# Probabilistic Multi-Criteria Decision-Making for Circularity Performance of Modern Methods of Construction Products

**Yiping Meng[a,b*], Sergio Cavalaro[a,**], Frozan Dizaye[a] Mohamed Osmani[a]**

[a] School of Architecture, Building & Civil Engineering, Loughborough University, Leicestershire, LE11 3TU, United Kingdom

[b] School of Computing, Engineering & Digital Technologies, Teesside University, Middlesbrough, Tees Valley, TS1 3BX, United Kingdom

* Correspondence to: Y. Meng, School of Architecture, Building & Civil Engineering, Loughborough University, Leicestershire, LE11 3TU, United Kingdom

** Correspondence to: S. Cavalaro, School of Architecture, Building & Civil Engineering, Loughborough University, Leicestershire, LE11 3TU, United Kingdom

E-mail addresses: y.meng@tees.ac.uk[1] (Y. Meng), s.cavalaro@lboro.ac.uk (S. Cavalaro), f.dizaye@lboro.ac.uk (F. Dizaye), m.osmani@lboro.ac.uk (M. Osmani)

[1] Present address: y.meng@tees.ac.uk (Y. Meng)

**Acknowledgements:**

This research is funded by EPSRC through the Interdisciplinary Circular Economy Centre for Mineral-Based Construction Materials from the UK Research and Innovation (EPSRC Reference: EP/V011820/1).

**Declaration of Interest statement:**

The authors declare that they have no known competing financial interests or personal relationships that could have appeared to influence the work reported in this paper.

## Abstract:

The construction industry faces increasingly more significant pressure to reduce resource consumption, minimise waste, and enhance environmental performance. Towards the transition to a circular economy in the construction industry, one of the challenges is the lack of a standardised assessment framework and methods to measure circularity at the product level. To support a more sustainable and circular construction industry through robust and enhanced scenario analysis, this paper integrates probabilistic analysis into the coupled assessment framework; this research addresses uncertainties associated with multiple criteria and diverse stakeholders in the construction industry to enable more robust decision-making support on both circularity and sustainability performance. By demonstrating the application in three real-world MMC products, the proposed framework offers a novel approach to simultaneously assess the circularity and sustainability of MMC products with robustness and objectiveness.

## 1. Introduction

The construction industry is one of the world's largest energy and raw materials consumers, responsible for around 40% of CO2 emissions [1] and nearly a third of all waste in the EU [2]. Hence, the construction industry is crucial in tackling climate change issues and reducing carbon emissions [3]. The concept of the circular economy refers to a new economic model which aims to reduce waste and pollution by improving efficiency while maintaining products and materials within the economy [4]. According to the Ellen MacArthur Foundation, the principles of CE are defined as eliminating waste and pollution, circulating products and materials at high value and regenerating nature [5]. These principles are often referred to as the 10 R's[1], which encompass a wide range of strategies to achieve circularity [6]. One of the aims of applying CE in construction is to minimise material usage and reduce waste at the construction and deconstruction stages [7], [8]. Through the adoption of circular economy principles, the construction industry can help achieve the United Nations' Sustainable Development Goals (SDGs), to drive environmentally friendly practices while creating positive economic and social impacts [9].

Modern Methods of Construction (MMC), which involves moving key construction processes to offsite locations, allow for greater control and precision in the construction process [10]. MMC is not restricted to the use of particular materials, allowing for greater flexibility in material selection, including steel, concrete, timber, and composites, among others, as opposed to traditional construction methods. As a result, MMC is considered to be a viable solution to mitigate the environmental impact of construction and increase productivity in the construction industry. [11], [12]. Empirical studies have demonstrated that incorporating MMC into construction processes can yield noteworthy reductions in embodied carbon (EC) and embodied energy (EE) associated with construction materials, as well as lower cumulative energy demand (CED), global warming potential (GWP), and construction waste. Additionally, MMC adoption has been shown to lead to shortened delivery times and increased labour productivity. Moreover, MMC has the potential to enhance material circularity by increasing the percentage of materials that can be reused or recycled [11].

Measuring the CE in the construction industry can help to stimulate and monitor the transition progress from a linear construction mode to a circular one, allowing for examining the effectiveness of the CE practices as well as highlighting the areas that require improvements [13], [14]. Evaluating the CE level can promote the application of CE's economic aspects in the construction sector by refining its strategies and approaches [15]. The main challenge of CE assessment is the lack of standardised metrics and indicators to measure circularity [16]. Current circularity assessment frameworks predominantly follow the CE principles at three distinct levels: macro, meso, and micro [17]. Nevertheless, there is still a need for the development of more refined frameworks and comprehensive research efforts that concentrate on CE assessment tools at design and product levels. Additionally, the nexus between sustainability and circularity remains undefined [18], [19], primarily due to the limited implementation of circular practices within the industries, especially for construction [20]. This ambiguity can result in a lack of accountability and transparency, making it challenging to assess the effectiveness of CE practices in minimising waste and maximising the material efficiency [21]. The most widely applied multi-criteria assessment for sustainability and/or circularity are developed based on the limited participants for questionnaires or interviews, which lacks objectivity and robustness [22].

[1] 10 Rs: Recover, Recycle, Repurpose, Remanufacture, Refurbish, Repair, Re-use, Reduce, Rethink and Refuse

To ultimately support a more sustainable and circular construction industry through robust and enhanced scenario analysis, this paper integrates probabilistic analysis into the assessment framework to compare the circularity performance with sustainability. This research aims to address uncertainties associated with multiple criteria and diverse stakeholders in the construction industry to enable more robust decision-making support.

## 2. Literature Review

### 2.1 Related work on circularity and sustainability relationship

Efforts to assess circularity performance have gained momentum globally, including established standards (AFNOR, 2018; UNI, 2022; ISO, 2024), government policies (European Commission, 2018; CLC, 2021; Eurostat, 2023) and organisational initiatives (EMF, 2015; Circle Economy, 2022; UKGBC, 2023). However, transitioning to a CE often entails trade-offs among environmental, economic, and social dimensions of sustainability, posing potential burdens rather than holistic benefits (Saidani *et al.*, 2024). As such, linking CE with sustainability performance in assessment methodologies is crucial for minimising burden-shifting impacts (Shevchenko *et al.*, 2024).

Several scholars have investigated the relationship between circularity and sustainability assessments, revealing both opportunities and challenges. For instance, Kravchenko, Pigosso and McAloone (2019) consolidated a database of sustainability indicators for CE strategies in manufacturing, followed by a dynamic selection process for aligning indicators with specific CE strategies and corporate contexts (Kravchenko, Pigosso and McAloone, 2020). However, the authors reliance on established frameworks (e.g. potting et al) limits adaptability to novel CE scenarios. Additionally, (Martinho, 2021) mapped CI against sustainability dimensions, highlighting conflicting metrics, a lack of uniformisation and a disproportionate focus on business and cities while neglecting water and energy aspects. In response, (Saidani and Kim, 2022), proposed a framework formalising connections between circularity and sustainability as beneficial, conditional, or scenario-dependent trade-offs. Similarly, (de Oliveira and Oliveira, 2023) analysed CE indicators across macro, meso, micro, and nano levels in relation to sustainability dimensions and lifecycle stages (take, make, use, recover), observing a predominant focus on material recycling and a lack of robustness to assess the sustainability of circular practices.

Recent developments include the *Circular Lifecycle Sustainability Assessment* framework by (Luthin, Traverso and Crawford, 2023, 2024), which integrates circularity with lifecycle assessment, lifecycle costing, and social LCA. This framework exposes trade-offs across sustainability dimensions. However, its reliance on metrics such as the Material Circularity Indicator (MCI) by the (EMF, 2015), limits its ability to link circularity with cradle-to-cradle lifecycle aspects such as (raw materials' supply risk) and its applicability to dynamic and sector-specific contexts (Hackenhaar *et al.*, 2024).

### 2.2 Related work on circularity and sustainability assessment for construction

For the construction sector, the most important question is what to measure in developing or selecting the circularity assessment framework. Various approaches have been developed to assess the performance of circularity within the construction sector, among which life cycle assessment (LCA) has emerged as a prominent method [23]. LCA can offer in-depth results on resource utilisation, waste generation, and environmental impact of the products along the CE process, providing opportunities for further targeted enhancement [24]. However, as LCA is developed referring to the linear economy mode, there are some deviations between the main concept of LCA-oriented eco-efficiency design and the environmental sustainability of CE, making LCA a less appropriate option to assess the CE in construction sector [25]. Moreover, LCA and life cycle cost (LCC) are implemented to assess the circularity of construction products quantitatively [26]. The results turn out that LCA and LCC not covering the cradle-to-cradle life cycle, lacking the information of construction stage and end-of-life of the building [27]. Apart from LCA and LCC, there are some other commonly used qualitative and quantitative assessment tools like material flow analysis (MFA), cost–benefit analysis (CBA) [13], [28]. Researchers have highlighted the need to determine the KPIs (key performance indicators) to have a comprehensive understanding for circularity assessment from the beginning of product design [28] to support different stakeholders of construction to make decisions on different design options [13].

A plethora of circularity and sustainability assessment methods have been populated across different CE levels: macro, meso, micro, nano, and cross-cutting summarized in Table 1.

Table 1 Circularity and sustainability assessment methods

| Systemic Level | Reference/ Name | Purpose | Method |
|---|---|---|---|
| Macro (Country) | (Eurostat, 2023)/ EU Circular Economy Monitoring Tool | Track CE progress in the EU | 11 statistical indicators based on official Eurostat data |
| | (Fatimah *et al.*, 2020)/ Multidimensional smart circular sustainable waste management framework | Assess the maturity of Indonesia's waste management systems using industry 4.0-based sustainable CE approach to achieve SDGs | Questionnaire survey and observation (4 cities) |
| | (Alfaro Navarro and Andrés Martínez, 2024)/ CE Implementation Index | Measure CE implementation in the EU countries and its link to SDGs | Mathematical modelling |
| Meso (Industrial networks) | (Ramírez-Rodríguez, Ormazabal and Jaca, 2024)/ Industrial symbiosis sustainability assessment tool | Map industrial symbiosis lifecycle sustainability assessment methods under the CE pathway | Conceptual modelling |
| | (Borbon-Galvez *et al.*, 2021)/ Cross-border construction waste management to | Assess sustainability of cross-border management of aggregates and construction and demolition waste between Italy and Switzerland | Semi-structured interviews |
| | (Cagno *et al.*, 2023)/ Multi-level performance measurement system | A multi-level performance measurement of sustainability, CE and industrial symbiosis | LCA |
| | (Fraccascia, Giannoccaro and Albino, 2021)/ Ecosystem-based industrial symbiosis indicators | Ecosystem-based industrial symbiosis performance indicators to measure resource efficiency and benefits | Enterprise Input-Output |
| Micro (Building) | (European Commission, 2018)/ Level(s) | Standardise sustainability reporting across a building's lifecycle | Reporting templates and sustainability indicators |
| | (Bronsvoort, 2021)/ Madaster Circularity Indicator | Assess building circularity based on material data of both technical and biological lifecycles | MFA |
| | (Han, Kalantari and Rajabifard, 2024)/ Building Information Modelling-Based Demolition Waste Management | Guide sustainability-oriented demolition waste decisions | LCA, Hybrid multi-criteria decision-aiding |
| | (Lederer and Blasenbauer, 2024)/ MFA-based sustainability assessment | Assess material flows and sustainability of CE scenarios for urban building stock | MFA |
| | (Khadim *et al.*, 2023)/ Whole-Building Circularity Indicator | Measure building circularity and sustainability | Mathematical modelling based on MCI (LCA complementary) |
| Nano (Products, | (Antwi-Afari *et al.*, 2023)/ Predictive Systemic Circularity | Evaluate circularity potential for modular steel slabs | LCA |
| materials) | Indicator | | |
| | (Jayawardana *et al.*, 2023)/ Modular Construction CE Assessment | Assess reuse and recycling potential of modular construction in developing economies | LCA |
| | (Bracquené, Dewulf and Duflou, 2020)/ Product Circularity Indicator | Measure the circularity of complex product supply chains | Mathematical modelling (comparison with MCI) |
| | (Muñoz, Hosseini and Crawford, 2024)/ 9R Circularity Index | Measure CE performance based on the 9R framework | MFA, AHP for weighted scoring |
| | (Steenmeijer *et al.*, 2024)/ Circularity Indicators for Recycling | Assess and classify recycling options for high-quality secondary resource flows | MFA |
| Cross-cutting | (ISO, 2024)/ ISO/ DIS 59020 | Provide generic, standardised methodology for measuring circularity performance and sustainability impacts | Circularity indicators complemented by LCA, MFA, Social Responsibility, and Life Cycle Costing standards. |
| | (Kazmi and Chakraborty, 2023)/ Indicators for Implementing Circularity | Identify and validate 78 circularity indicators for construction | Delphi survey (30 experts) |
| | (Platform CB'23, 2022)/ Platform CB'23 - circularity measurement framework | Guide Dutch construction sector's circularity measurements | Stakeholder-based qualitative framework |

## 3. Probabilistic MCDM for Coupled Circularity and Sustainability

### 3.1 Four-layer Framework Establishment

Drawing from a comprehensive literature review of multiple assessment tools and methods sustainability and/or circularity, the proposed assessment framework will incorporate both qualitative and quantitative factors. Referring to the widely recognised sustainability assessment frameworks and main categories of CE at product level [29], this study will include economic, social, and environmental aspects as the factors at the requirement level ($B_i$), by selecting the overlapped part of sustainability and circularity.

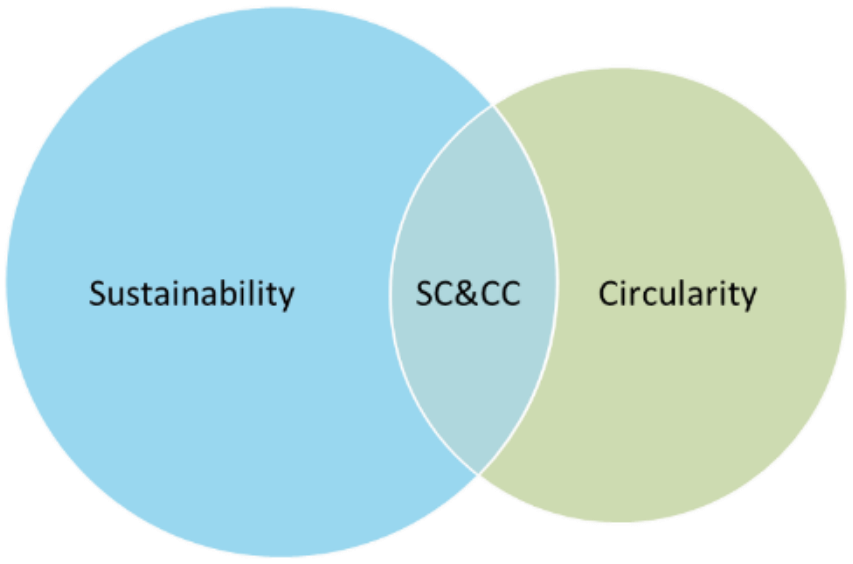

Figure 1 Venn Diagram for the selected criteria boundary

Recognising the unique attributes of MMC, technology is incorporated as an additional, the fourth factor at this level. The factors at the criteria level ($C_j$) are derived from related scholarly works, resulting in a total of 12 criteria, with simplicity and comprehensiveness. The factors at the criteria level ($C_{jk}$) are informed by LCC, the cradle-to-cradle framework, principles of CE, and the very definition of MMC. Decision-makers, however, retain the flexibility to further refine this framework to align with their specific objectives - they can either introduce additional criteria or pare down existing ones. This flexibility ensures the proposed framework's broader adaptability and relevance across a variety of MMC contexts.

Table 1 Primary List of Criteria and Indicators

| | Requirement level | Criteria Level | Indicator Level | Impact factors/explanations | References |
|---|---|---|---|---|---|
| | Economics ($B_1$) | Costs ($C_1$) | Pre-construction costs $C_{11}$ | coordination costs | [30] |
| | | | Initial Production and assembly cost $C_{12}$ | material costs | [31], [32], [33], [34] |
| | | | | transportation costs | |
| | | | | Labour costs | |
| | | | | equipment costs | |
| | | | Maintenance cost $C_{13}$ | Cost for repair the equipment | [35], [36] |
| | | | End-of-life costs $C_{14}$ | disposal costs | [37] |
| | | | | demolition costs | |
| | | | | waste treatment costs | [38] |
| | | | | recycle costs | |
| | | Time $C_2$ | lead time $C_{21}$ | | [37] |
| $A$ | | | Initial Production and assembly time $C_{22}$ | | [39] |
| | | | Time for end-of-life process $C_{23}$ | Deconstruction time | [40] |
| | | | | Disposal time | |
| | | | | recycle time | |
| | | Investment Return $C_3$ | the speed of return on investment $C_{31}$ | the speed of return on investment | [37] |
| | Social ($B_2$) | Health and Safety $C_4$ | Users' health $C_{41}$ | indoor quality | [37] |
| | | | Labours' safety $C_{42}$ | working condition | [32], [37] |
| | | | | Labour risk of accidents during construction | [32], [37] |
| | | Surrounding impact $C_5$ | onsite disturbance $C_{51}$ | duration of time | [36], [37] |
| | | | Service change $C_{52}$ | service quality | [37] |
| | | | | service capacity | [41] |
| | Environment ($B_3$) | Material consumption $C_6$ | | | [39] |

| | | | | |
|---|---|---|---|---|
| | Energy consumption $C_7$ | | | [37], [42] |
| | Waste Generation $C_8$ | | | [37], [39] |
| | CO2 emission $C_9$ | | | [43] |
| Technology ($B_4$) | Product adaptability $C_{10}$ | Modularisation level $C_{101}$ | | [44] |
| | | Flexibility $C_{102}$ | | [44] |
| | | recycling and demolition ability $C_{103}$ | disassembly | |
| | | | Recyclable/renewable contents | [37] |
| | | | Recyclable/reusable elements | [32] |
| | Technical quality $C_{11}$ | Durability $C_{111}$ | | [37], [45] |
| | | usage efficiency $C_{112}$ | | [38] |
| | | Defects and damages $C_{113}$ | | [37], [46] |
| | | Accuracy $C_{114}$ | | [47], [48] |
| | Technical capacity $C_{12}$ | prefabrication degree $C_{121}$ | | [49], [50], [51] |
| | | The scale of mass production $C_{122}$ | | [34] |
| | | Productivity $C_{123}$ | | [52] |

3.2 Questionnaire Design

Building upon the established four-level assessment framework, an interactive online questionnaire survey was designed to probe the importance level of the of 12 distinct criteria. The main aim of the survey is to understand the deviation between sustainability and circularity paradigm from the perspective of diverse stakeholders within the MMC domain. For each criterion, respondents were prompted to ascribe an importance rating on a scale of 0 to 10, where 0 corresponded to 'not relevant' and 10 to 'most important'. The other ratings (1-2, 3-4, 5-6, 7-8) denoted varying degrees of importance ranging from 'least important' to 'fairly important'. To facilitate comprehension, each criterion was accompanied by a brief description outlining the involved factors and their definitions. Notably, respondents were also invited to propose any additional criteria they considered important but unlisted, thereby affording the study the benefit of their unique insights and expertise. The questionnaire was disseminated during meetings and workshops organised by the Circular Economy Center, thereby ensuring a relevant and informed pool of respondents.

3.3 Probabilistic Method

In this research, the weights of the third layer – **Criteria level** is calculated by AHP based on the data collected from the questionnaires. For the fourth layer-the indicator level, 85 samples from six group of stakeholders could not provide data with strong validity. To overcome uncertainty due to varying expert judgments, differing perspectives of stakeholders, and incomplete information, the Monte Carlo simulation is utilised to generate the weights for the fourth layer. Monte Carlo is based on the repetitive computational execution of a certain number of random deterministic scenarios to model the probability of different outcomes [53]. Sampling method is often used together to reduce the repeatability of random sampling process of Monte Carlo [54]. One of the sampling improvements for Monte Carlo is Latin Hypercube Sampling (LHS), which can simplify the simulation with high accuracy [55]. The steps for generating the weights for the fourth layer using LSH-Monte Carlo follows:

(1) Determine the value range. Define the value for the weights of the indicators under each criteria. The value is following the constraints:

a) s.t.

b) $$\begin{cases} \sum_{j=1}^{m} w_{c_{ij}} = 1 \\ 0.1 \leq min_{w_{c_{ij}}} < \frac{1}{m} \\ \frac{1}{m} \leq max_{w_{c_{ij}}} \leq 1.1 - 0.1 \times m \\ m \in \{2,3,4,5\} \end{cases} \quad (1.10)$$

c) Where $i = 1,2,\ldots,12$, and $j = 1,2,\ldots m$ and $m$ is the number of indicators in criteria $c_i$.

d) The minimum value of $w_{c_{ij}}$ is restricted to 0.1, to avoid the situation that the weights for some indicators are nearly zero. Then the maximum value for the $w_{c_{ij}}$ can also be limited to make sure the sum of weights $w_{c_{i1}}, w_{c_{i2}}, \ldots, w_{c_{im}}$ of criteria $c_i$ equals to one.

(2) Determine the distribution of the probability. This step is operated on each indicator. After defining the range, the distribution for the weight value is determined to generate the value. In this case, we want 1000 sets of weights for the indicators. Then the distribution for the value is divided into 1000 equal probability intervals.

(3) Random sampling process. This step is operated within one indicator. Based on the 1000 segmented probability intervals, one value is selected randomly from each interval. For the $m$ indicators in the criteria $c_i$, this process would operate $m$ times to get $w_{c_{i1}}, w_{c_{i2}}, \ldots, w_{c_{im}}$.

(4) Normalisation. This step is to make sure the sum of randomly selected weights $w_{c_{i1}}, w_{c_{i2}}, \ldots, w_{c_{im}}$ is $1$. The normalisation will generate the normalised weights $w'_{c_{i1}}, w'_{c_{i2}}, \ldots, w'_{c_{im}}$ (equation 1.11).

(5) $$w'_{c_{ij}} = \frac{w_{c_{ij}}}{\sum_{j=1}^{m} w_{c_{ij}}} \quad (1.11)$$

(6) Construct the weights. After the normalisation process, one set of weights is generated $[w'_{c_{i1}} \quad w'_{c_{i2}} \quad \cdots \quad w'_{c_{im}}]_{1\times m}$.

(7) Replicate the sampling process. Repeat the step c-e 1000 times to get the weight vector for the $m$ indicators of criteria $c_i$. The row of the matrix represents the $m$ weights for the indicators and the column represents the $1000$ selected weights for $c_{ij}$.

a) $\begin{bmatrix} w_{c_{i1}}^{1} & \cdots & w_{c_{im}}^{1} \\ \vdots & \ddots & \vdots \\ w_{c_{i1}}^{1000} & \cdots & w_{c_{im}}^{1000} \end{bmatrix}_{1000\times m}$

(8) Construct the weight vectors for the fourth layer. Repeat the step c-e $n$ times to generate all the indicators of the 12 criteria. After this step a final weight matrix is obtained. The size for the matrix is $1000 \times \sum_{i=1}^{n} m_i$.

a) $\begin{bmatrix} w_{c_{11}}^{1} & \cdots & w_{c_{12m}}^{1} \\ \vdots & \ddots & \vdots \\ w_{c_{11}}^{1000} & \cdots & w_{12m}^{1000} \end{bmatrix}_{1000\times\sum_{i=1}^{n} m_i}$

(9) Finalise the weight matrix for sustainability and circularity. Based on the hierarchy decision tree, finalise the 1000 sets of weights for sustainability and circularity based on the number of criteria and indicators.

4. Results for Tree-diagram Determination

The questionnaire was shared in the events like workshop and hackathon organised by Circular Economy Center. The participants were stakeholders from different areas of sustainability and construction. From December 1st 2022 to January 31st 2023, 85 questionnaires from six groups of stakeholders were collected. And the distribution of the participants is shown in Figure 2.

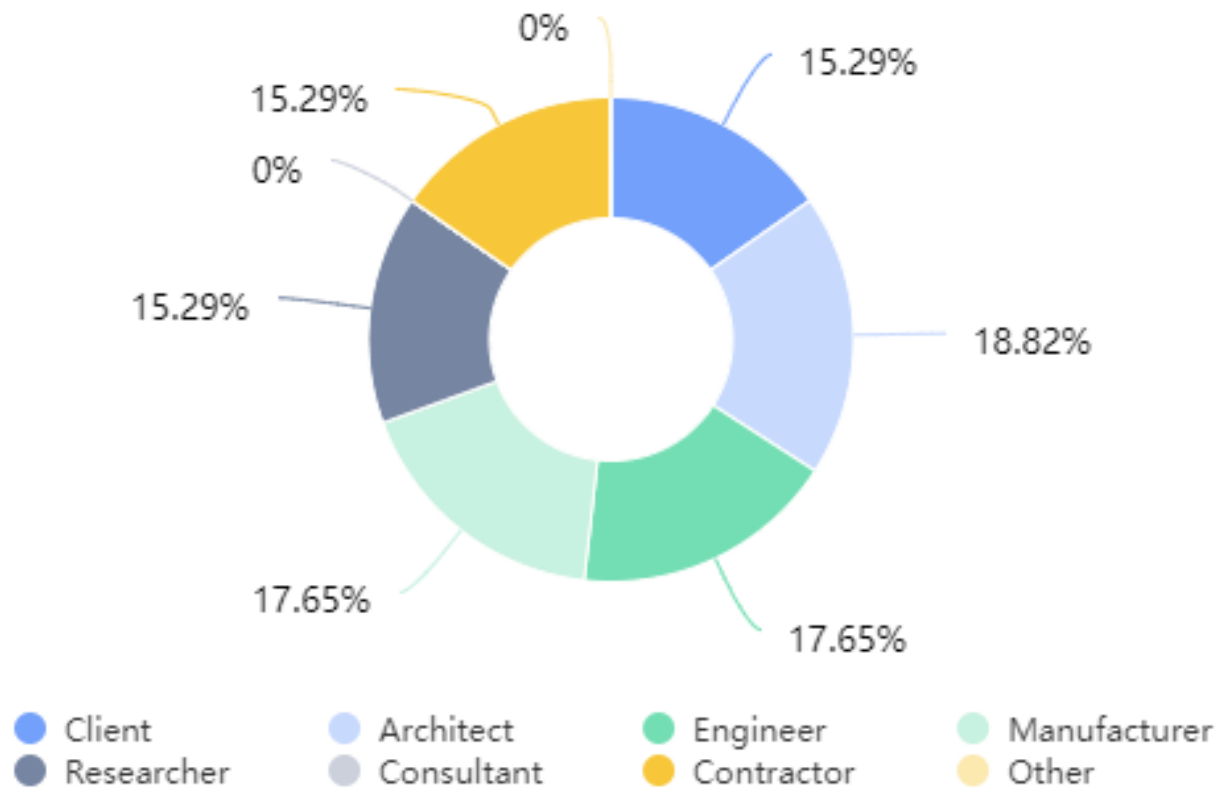

Figure 2 Work categories of survey participants

Based on the decision tree and the collected data from the questionnaire, the AHP is operated to calculate the weights for the 12 criteria. Table 2 and Table 3 depict the calculated weights for the 12 criteria for each of the six identified stakeholder groups involved in this study. It's essential to note that the column labelled 'General' in these tables represents the combined weights from all 85 participant responses. This aggregate data provides an overall picture of the perceived importance of each criterion across all stakeholders, offering a more holistic view of the assessment framework's applicability.Table 2 AHP Weights for different stakeholders for

Sustainability

| | Stakeholders | | | | | | |
|---|---|---|---|---|---|---|---|
| Criteria | Architect | Engineer | Manufacture | Contractor | Researcher | Client | General |
| Cost | 35.650% | 35.241% | 34.392% | 34.821% | 34.375% | 34.018% | 34.740% |
| Time | 32.024% | 33.735% | 33.862% | 36.310% | 32.188% | 33.724% | 33.660% |
| Investment Return | 32.326% | 31.024% | 31.746% | 28.869% | 33.437% | 32.258% | 31.600% |
| Health and Safety | 49.774% | 50.463% | 51.628% | 48.276% | 51.415% | 49.758% | 50.192% |
| Surrounding Impact | 50.226% | 49.537% | 48.372% | 51.724% | 48.585% | 50.242% | 49.808% |
| Material | 23.13% | 33.138% | 23.592% | 23.732% | 24.541% | 23.950% | 23.743% |
| Energy | 24.911% | 34.897% | 24.648% | 25.152% | 25.000% | 27.101% | 25.268% |
| Waste | 24.377% | 31.965% | 25.704% | 24.949% | 22.936% | 23.319% | 24.424% |
| CO2 | 27.580% | 33.138% | 26.056% | 26.166% | 27.523% | 25.630% | 26.565% |
| Adaptability | 35.342% | 33.138% | 32.471% | 31.464% | 32.692% | 36.486% | 33.585% |
| Technical quality | 31.507% | 34.897% | 34.195% | 35.514% | 34.295% | 33.446% | 33.938% |
| Technical Capacity | 33.151% | 31.965% | 33.333% | 33.022% | 33.013% | 30.068% | 32.476% |

Table 3 Weights for different stakeholders for Circularity

| | Stakeholders | | | | | | |
|---|---|---|---|---|---|---|---|
| Criteria | Architect | Engineer | Manufacture | Contractor | Researcher | Client | General |
| Cost | 34.043% | 35.018% | 35.599% | 35.000% | 34.268% | 34.539% | 34.725% |
| Time | 31.915% | 31.408% | 32.362% | 36.071% | 32.399% | 30.921% | 32.473% |
| Investment Return | 34.043% | 33.574% | 32.039% | 28.929% | 33.333% | 34.539% | 32.802% |
| Health and Safety | 51.596% | 55.276% | 54.737% | 47.596% | 51.456% | 51.323% | 51.949% |
| Surrounding Impact | 48.404% | 44.724% | 45.263% | 52.404% | 48.544% | 48.677% | 48.051% |
| Material | 27.091% | 25.475% | 26.326% | 25.482% | 25.624% | 24.886% | 25.864% |
| Energy | 22.182% | 23.194% | 21.780% | 24.197% | 23.810% | 25.799% | 23.390% |
| Waste | 28.727% | 27.757% | 29.735% | 27.623% | 24.263% | 26.256% | 27.525% |
| CO2 | 22.000% | 23.574% | 22.159% | 22.698% | 26.304% | 23.059% | 23.220% |
| Adaptability | 36.893% | 35.109% | 36.516% | 33.043% | 31.858% | 35.331% | 34.922% |
| Technical quality | 32.039% | 33.656% | 31.981% | 35.942% | 34.513% | 35.016% | 33.719% |
| Technical Capacity | 31.068% | 31.235% | 31.504% | 31.014% | 33.628% | 29.653% | 31.359% |

Based on the updated reliability analysis for the collected data, indicator waste generation is removed for the circularity tree diagram (details see supplementary material). The final decision tree for sustainability and circularity is adjusted, shown in Figure 3- Figure 4.

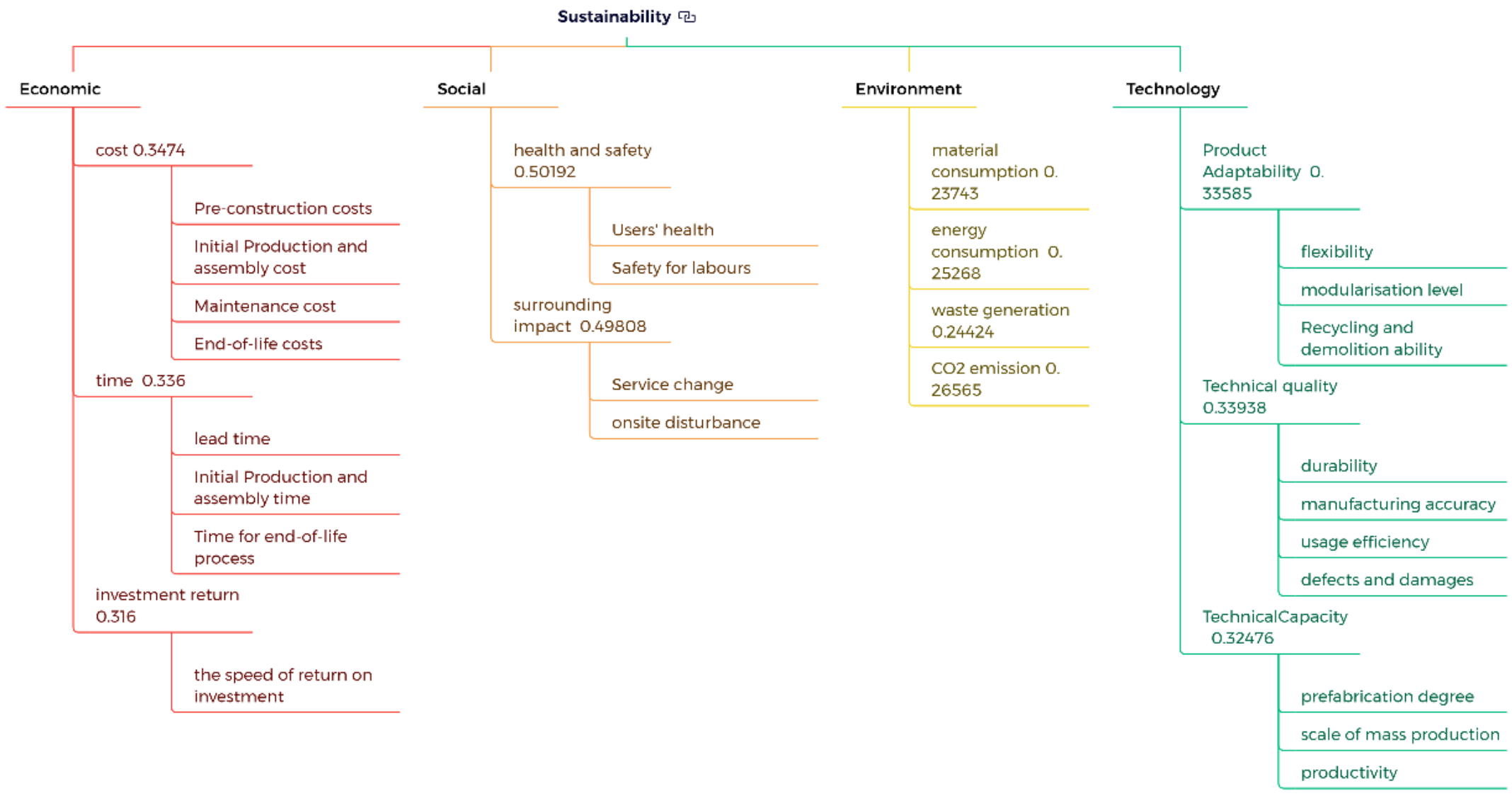

Figure 3 Finalise Decision Tree for Sustainability

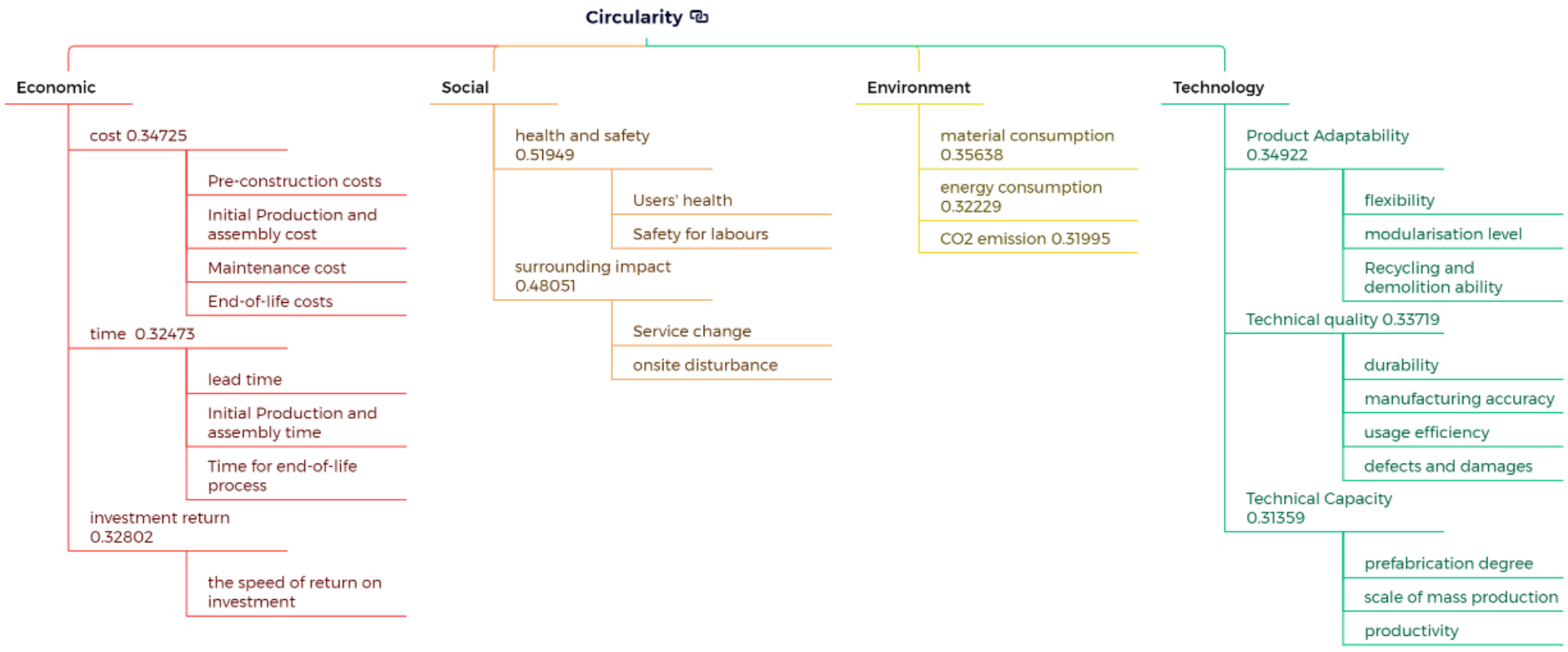

Figure 4 Finalise Decision Tree for Circularity

# 5. Application on Multi-criteria Decision Making

## 5.1 MIVES Application

Based on the weights obtained from AHP and simulation results of Monte Carlo, the framework is established with weights. To testify the coupled assessment framework, one multi-criteria decision-making (MCDM) method named MIVES (Modelo Integrado de Valor para una Evaluación Sostenible) [56] is utilised. MIVES has been applied in different assessment scenarios [57], [58] including for MMC in structural level [59]. MIVES can uniform quantitative and qualitative data into value between 0-1 using the equation (1.8)-(1.9). Equation (1.8) is the value function determined by parameter $F$, $C$ , $K$ and $B$. In this case, we select three MMC scenarios to apply the probabilistic assessment framework. As we generate 1000 sets of weights for the fourth layer, the ranking results of the three scenarios are value with probability. MIVES has defined eight types of curve shapes, namely convex, concave, linear, and S-shape. The type of curve shape depends on whether the curve shows an increasing or decreasing trend.

$$IV_i = B_i * \left[1 - e - Ki * \left(\frac{|X - X_{min_i}|}{C_i}\right)^{F_i}\right] (1.8)$$

$$B_i = \left[1 - e - Ki * \left(\frac{|X_{max_i} - X_{min_i}|}{C_i}\right)^{F_i}\right]^{-1} (1.9)$$

Three cases of MMC products in real world are selected and are defined as three scenarios for comparisons, namely:

Scenario 1 (S1) -A1 Panel from Spec Wall

Scenario 2 (S2) - HardiePanel®(Façade Panel)

Scenario (S3) - Precast Wall (DanElement)

These three MMC products are all made by precasting. These three products are in planar shape which can be applied as wall panel and façade.

Table 4 Value for MIVES

| Indicator | $X_{min}$ | $X_{max}$ | $F_i$ | $C_i$ | $K_i$ | Shape | S1 | S2 | S3 |
|---|---|---|---|---|---|---|---|---|---|
| $C_{11}$ | 0 | 10 | 1 | 6 | 1 | D-Convex | 6 | 4 | 4 |
| $C_{12}$ | 0 | 10 | 1 | 6 | 1 | D-Convex | 4 | 6 | 4 |
| $C_{13}$ | 0 | 10 | 1.2 | 4 | 1 | D-S | 6 | 6 | 8 |
| $C_{14}$ | 0 | 10 | 1.2 | 6 | 1 | D-S | 6 | 6 | 4 |
| $C_{21}$ | 0 | 10 | 1 | 4 | 1 | D-Convex | 2 | 6 | 4 |
| $C_{22}$ | 0 | 10 | 1 | 4 | 1 | D-Convex | 4 | 6 | 8 |
| $C_{23}$ | 0 | 10 | 1 | 6 | 1 | D-Convex | 8 | 6 | 4 |
| $C_{31}$ | 0 | 10 | 1.2 | 6 | 1 | D-S | 8 | 6 | 4 |
| $C_{41}$ | 0 | 10 | 1 | 8 | 1 | I-Convex | 8 | 4 | 6 |
| $C_{42}$ | 0 | 10 | 1 | 8 | 1 | I-Convex | 8 | 8 | 8 |
| $C_{51}$ | 0 | 10 | 1.2 | 4 | 1 | D-S | 4 | 4 | 6 |
| $C_{52}$ | 0 | 10 | 1 | 6 | 1 | I-Convex | 6 | 6 | 8 |
| $C_6$ | 0 | 10 | 1 | 4 | 1 | D-Convex | 6 | 4 | 4 |
| $C_7$ | 0 | 10 | 1 | 2 | 1 | D-Convex | 6 | 8 | 8 |
| $C_8$ | 0 | 10 | 3 | 4 | 0.01 | D-Concave | 6 | 4 | 4 |
| $C_9$ | 0 | 10 | 3 | 4 | 0.01 | D-Concave | 6 | 6 | 6 |
| $C_{101}$ | 0 | 10 | 1.2 | 4 | 1 | I-S | 8 | 8 | 6 |
| $C_{102}$ | 0 | 10 | 1.2 | 4 | 1 | I-S | 4 | 6 | 6 |
| $C_{103}$ | 0 | 10 | 1 | 6 | 1 | I-Convex | 4 | 6 | 6 |
| $C_{111}$ | 0 | 10 | 1 | 6 | 1 | I-Convex | 8 | 6 | 6 |
| $C_{112}$ | 0 | 10 | 1 | 6 | 1 | I-S | 8 | 6 | 4 |
| $C_{113}$ | 0 | 10 | 3 | 4 | 0.01 | D-Concave | 6 | 6 | 4 |
| $C_{114}$ | 0 | 10 | 1 | 8 | 1 | I-Convex | 6 | 8 | 6 |
| $C_{121}$ | 0 | 10 | 3 | 4 | 0.01 | I-Concave | 6 | 6 | 8 |
| $C_{122}$ | 0 | 10 | 1.2 | 4 | 1 | I-S | 8 | 6 | 6 |
| $C_{123}$ | 0 | 10 | 1 | 6 | 1 | I-Convex | 6 | 8 | 6 |

## 5.2 Probability Results for the Sustainability

We conducted 1000 experiments and obtained 1000 sustainability values for three different scenarios. The PDF and CDF for these values are shown in the figure, indicating the reliability of the simulation data. The value interval for S1, S2, and S3 were found to be 0.52-0.565, 0.525-0.57, and 0.53-0.58, respectively.

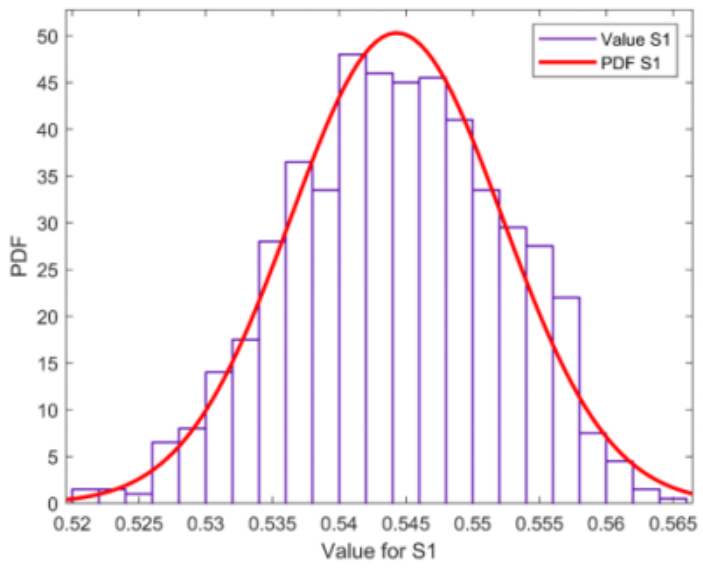

(a) PDF for the sustainability value of S1

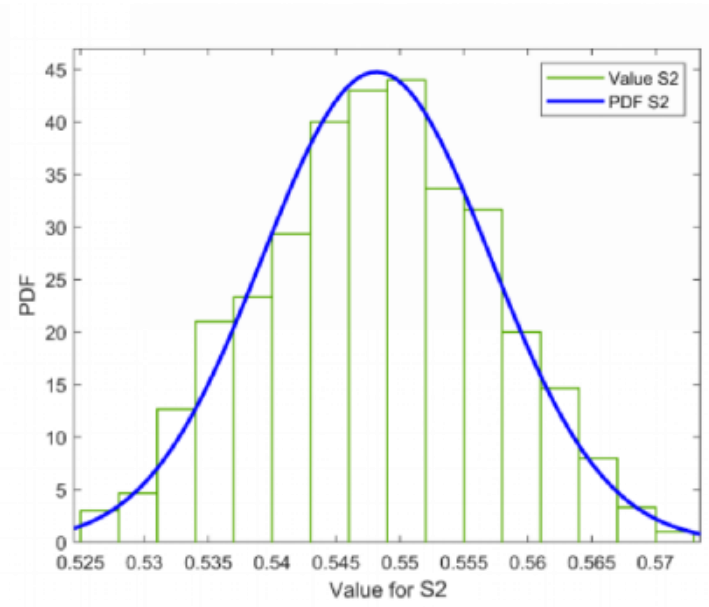

(b) PDF for the sustainability value of S2

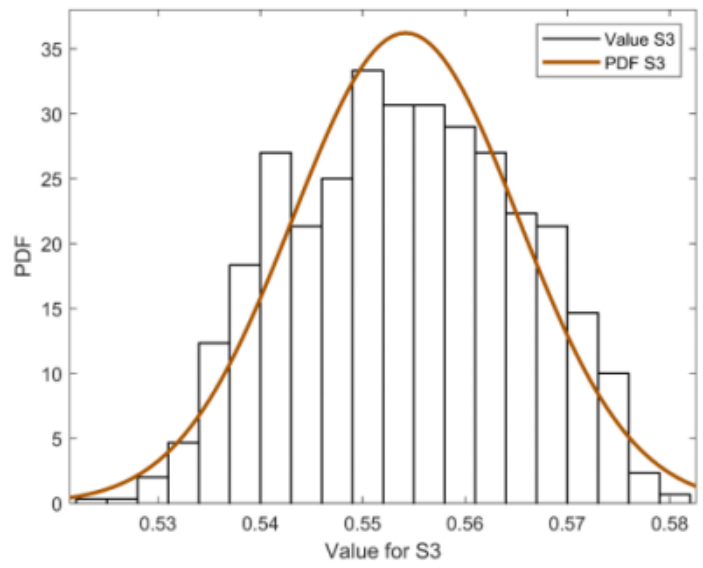

(c) PDF for the sustainability value of S3

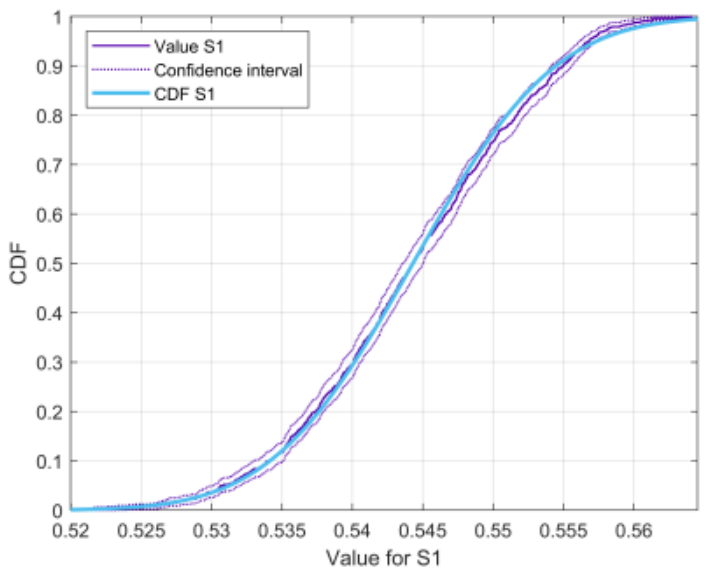

CDF for the sustainability value of S1

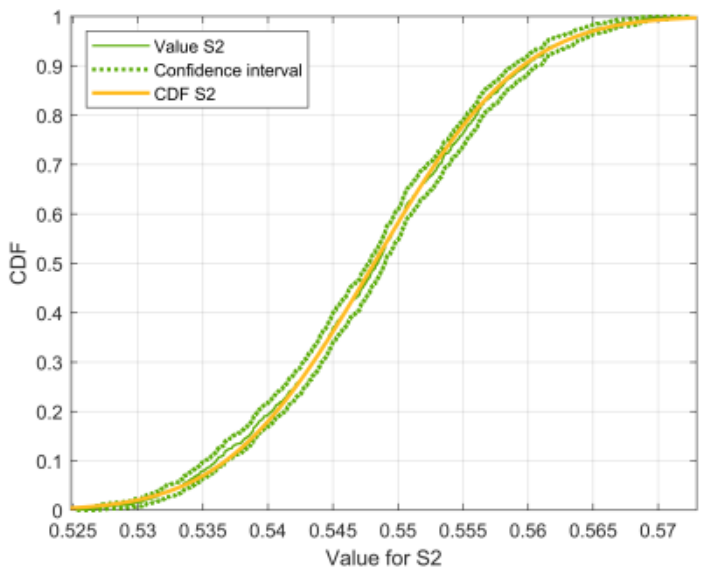

CDF for the sustainability value of S2

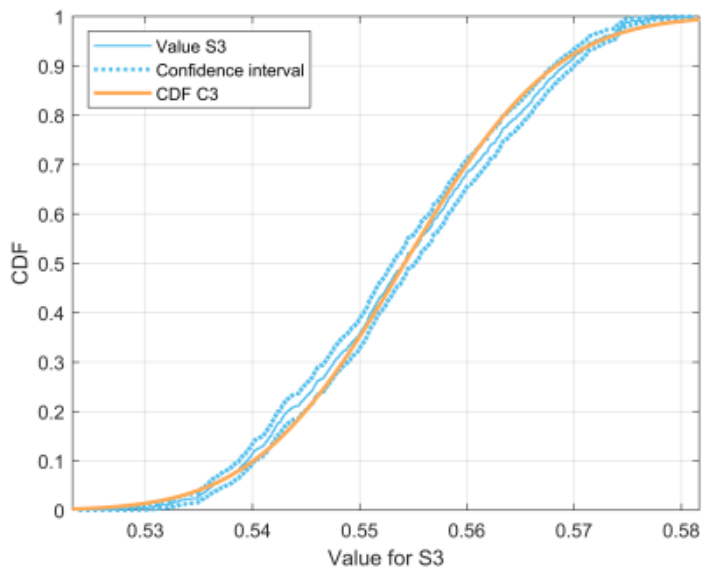

CDF for the sustainability value of S3

Figure 5 PDF and CDF for sustainability value for three scenarios

The mean values of three MMC product scenarios for sustainability level are summarised at first. According to the Figure 6, Scenario 3 has the highest sustainability level of 0.6346, followed by Scenario 2 with 0.6123. Scenario 1 has the lowest sustainability level with 0.5397.

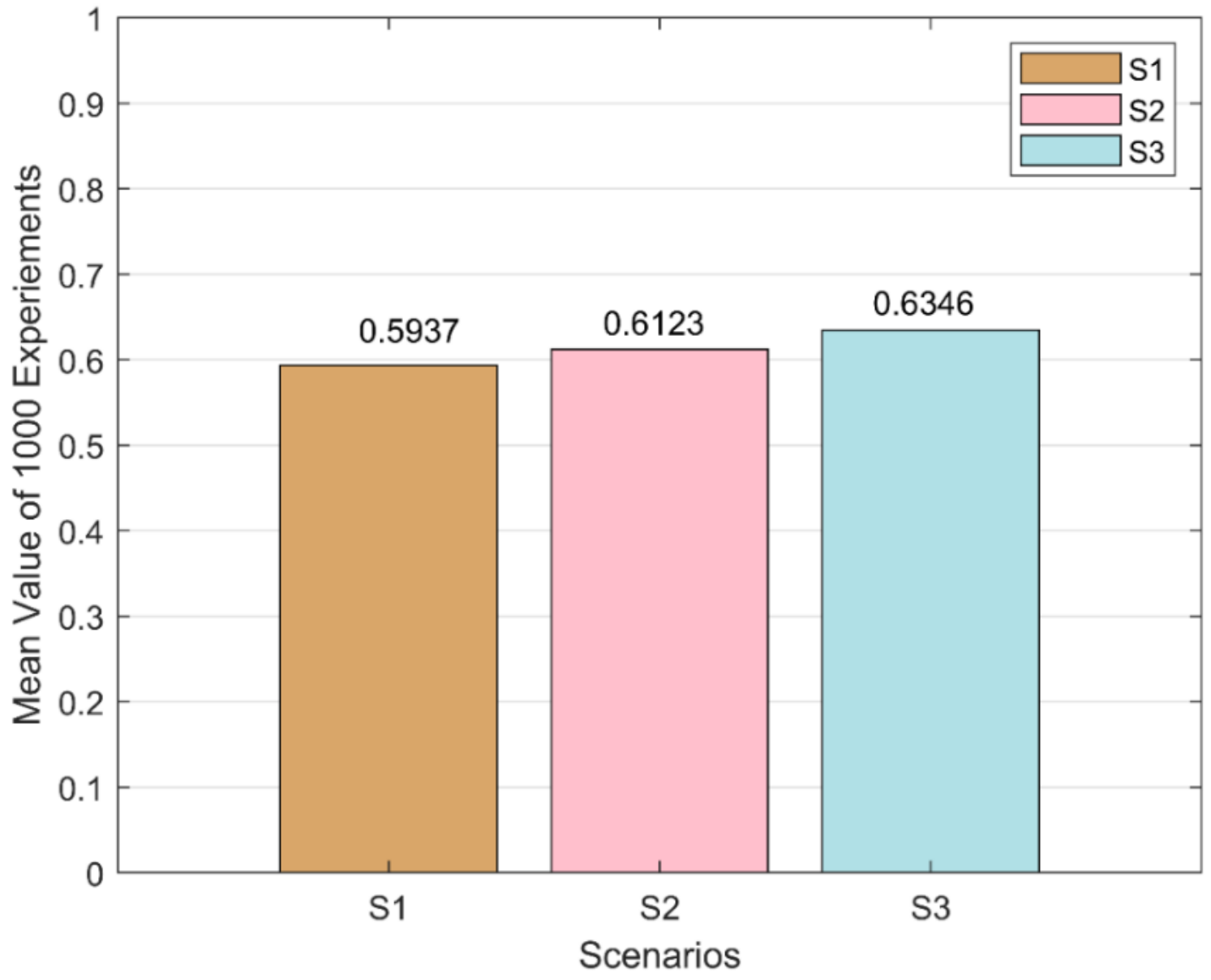

Figure 6 Mean Value of Sustainability for Three Scenarios

When analyzing sustainability performances across four different requirement levels, shown in Figure 7, S2 stands out in terms of environment and technological performance. S1, on the other hand, excels in social performance, while S3 performs best in the economic aspect. When considering the ratio of scores for each level in relation to the overall sustainability value, the social sector is the most significant contributor, followed by technology, economics, and the environment in descending order.

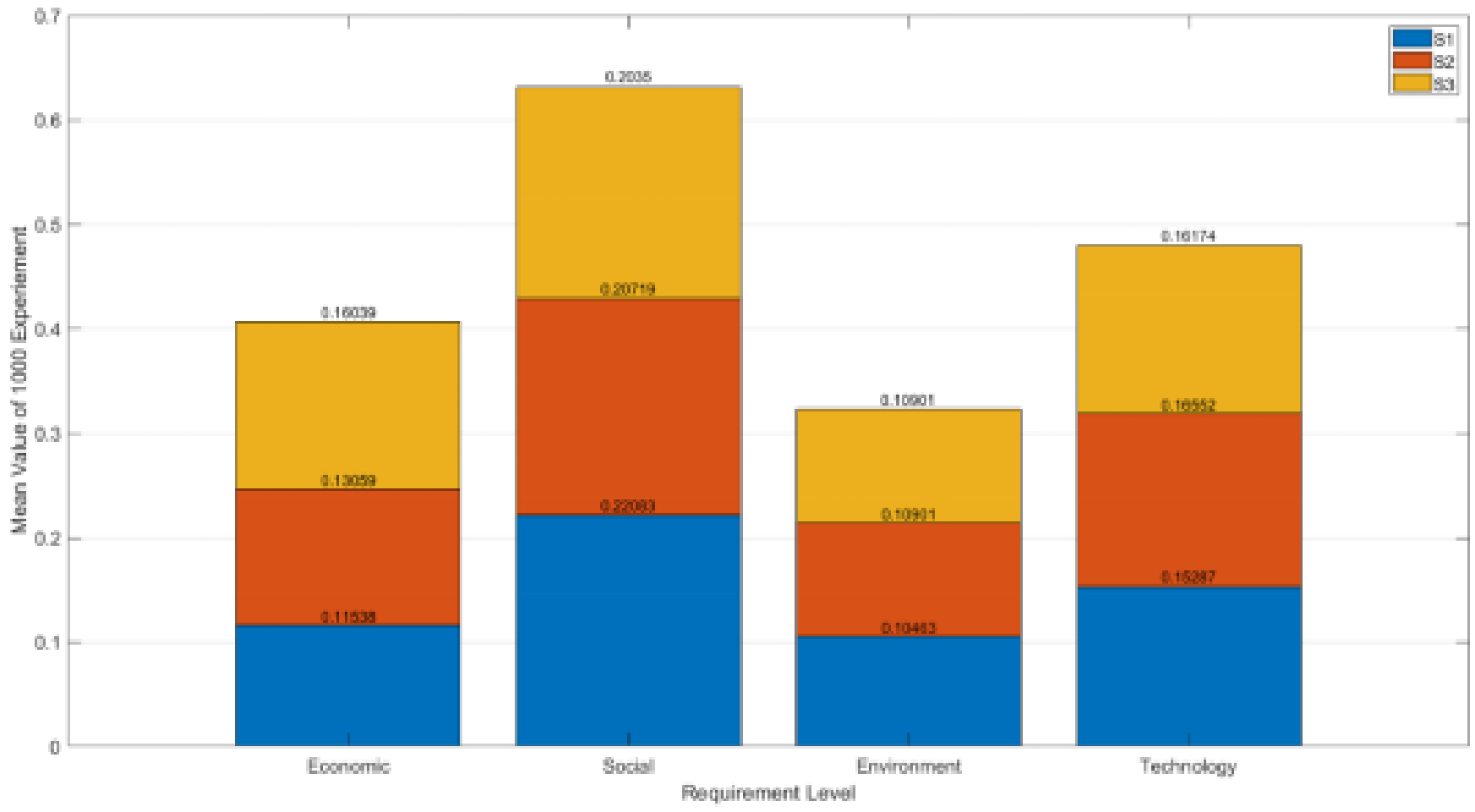

Figure 7  Stack bar for the mean value of sustainability of the three scenarios in four requirements

The Figure demonstrates a detailed analysis of 12 criteria. In terms of the economy, S1 performs the best in C1 and C2, whereas S3 has the best performance in C3. For the social aspect, S1 scores the highest in C4 criteria, and S1 and S2 have similar performances in C5. All three scenarios have similar environmental performances in C7 and C9. S1 has the best result for C6, and S3 has the best result in C8. In terms of technology, S3 has the highest value for C10 and C12, while S1 has the highest in C11.

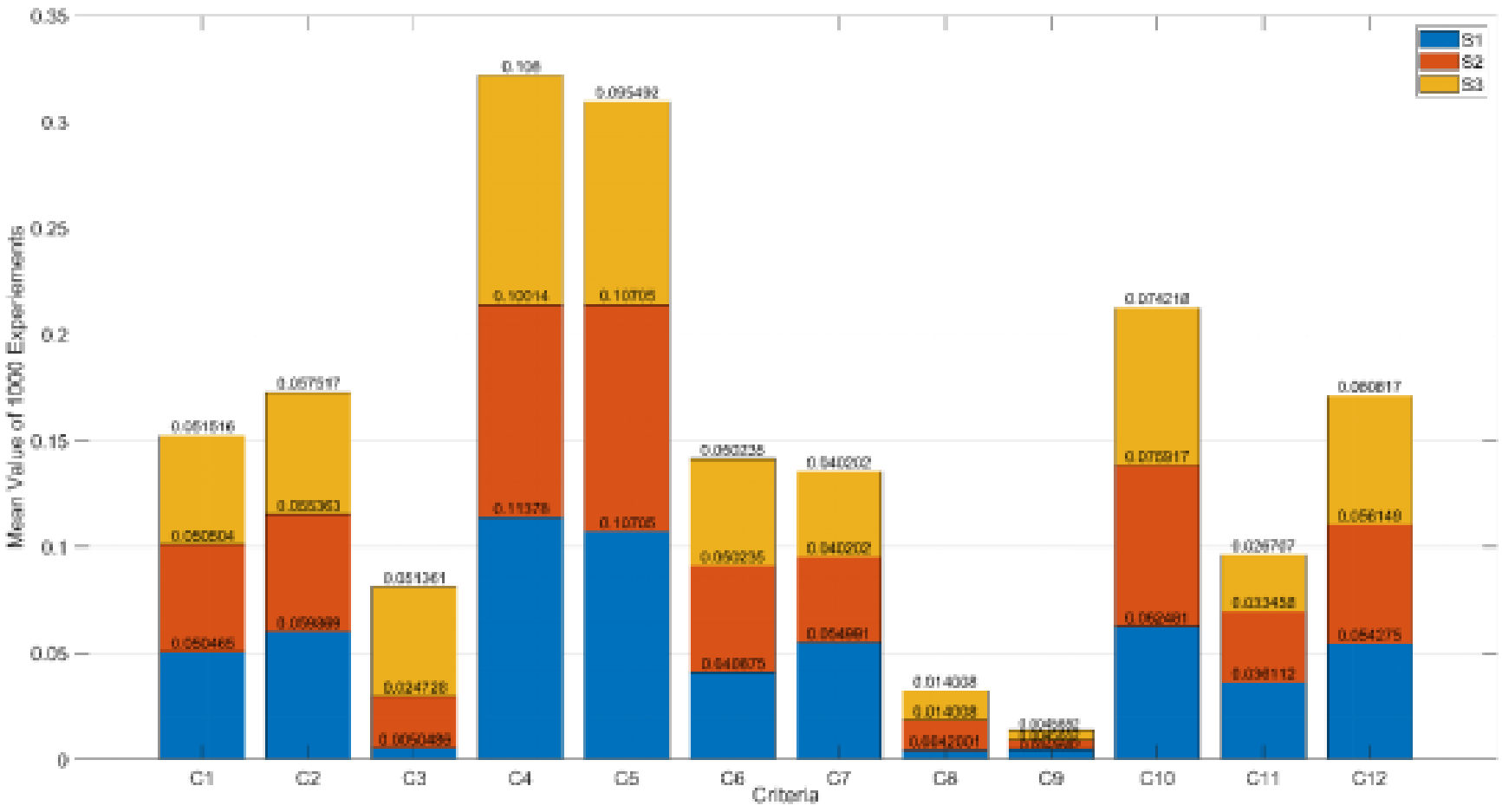

Figure 8 Stack bar for the mean value of sustainability comparison between 12 criteria

Upon analysis of the Figure 9 and Figure 10, it is evident that each of the three scenarios, S1, S2, and S3, excel in different areas. Specifically, S1 performs exceptionally well in the social category and C4 criteria, while S2 demonstrates notable success in the technology category and C5 criteria. Finally, S3 stands out in the economic category and C4 criteria. Overall, it is apparent that each scenario has its unique strengths and areas for improvement.

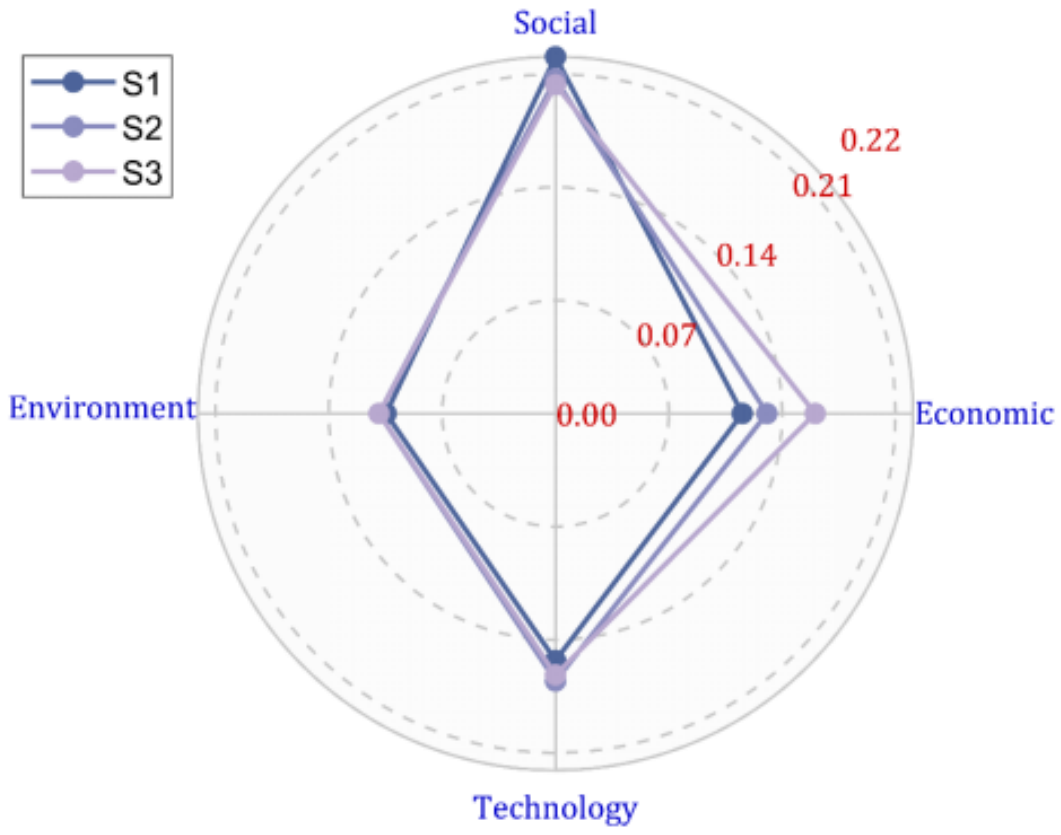

Figure 9 Radar chart for mean sustainability value of four requirements of three Scenarios

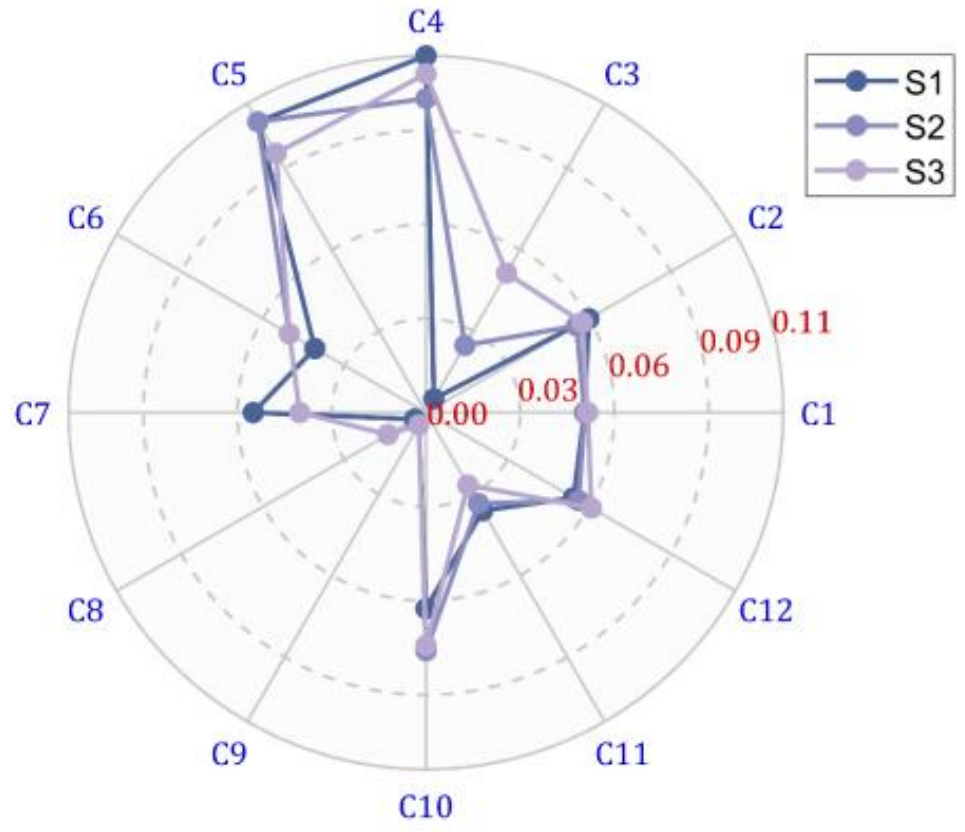

Figure 10 Radar chart for mean sustainability value of 12 criteria of three Scenarios

Compared to the mean value, Figure 11- Figure 12 show the probability of three scenarios rank as first and third based on their sustainability values. The detailed probability is calculated. Based on the results, it appears that S3 is most likely to rank first with a probability of 90.4%. On the other hand, S1 has the highest likelihood of ranking third with a probability of 93.9%, when taking into account the overall sustainability performance.

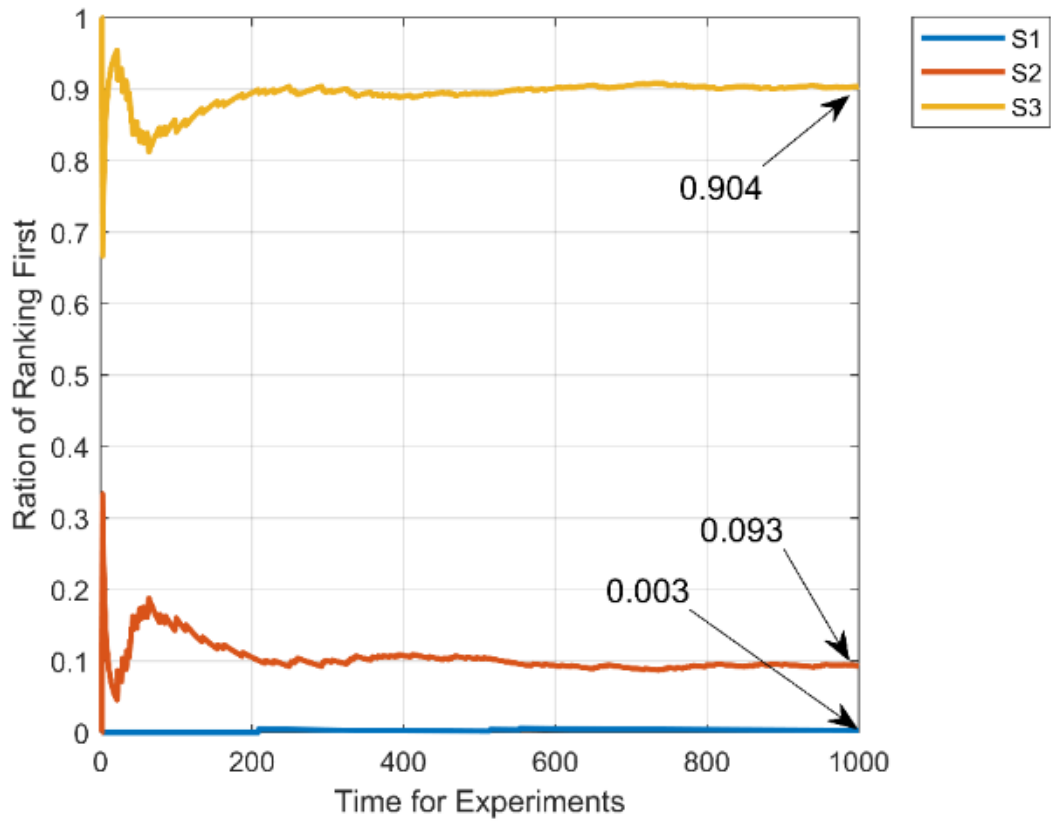

Figure 11 Ratio of ranking as first for three scenarios in sustainability (minimum weight 0.1)

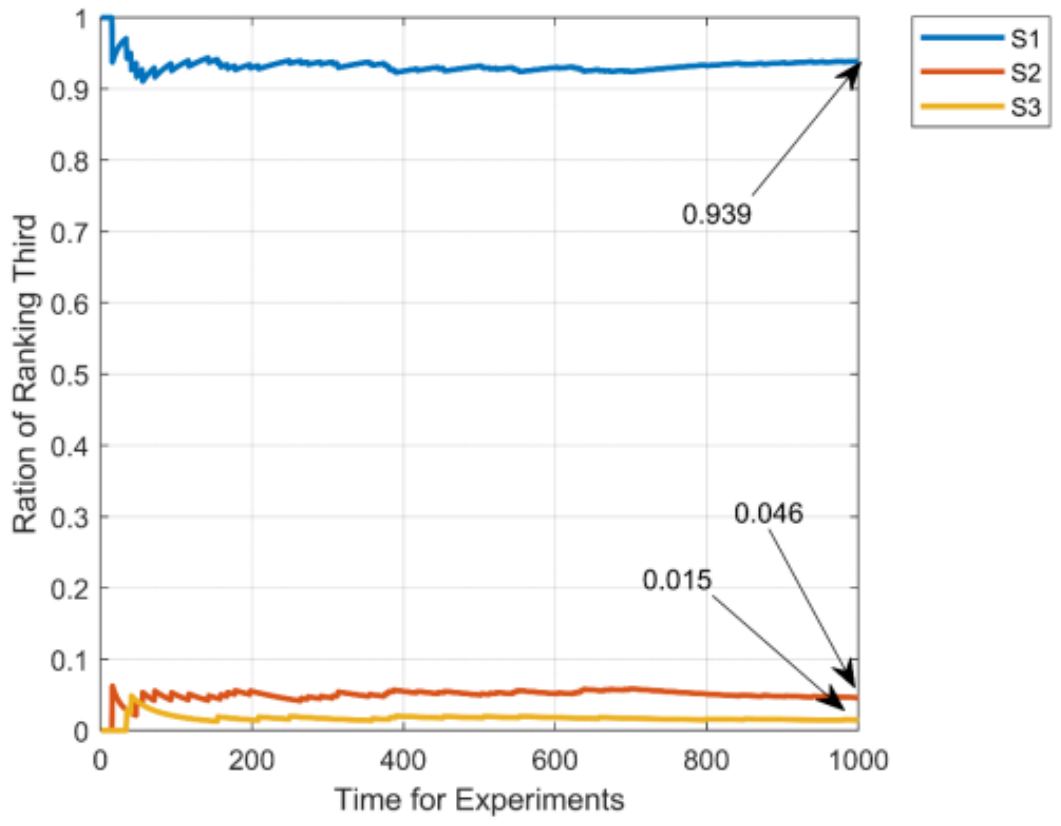

Figure 12 Ratio of ranking as third for three scenarios in Sustainability (minimum weight 0.1)

The figures below (Figure 13) demonstrate the probability of ranking first in four sustainable requirements. S3 has a 100% probability of ranking first in economic sustainability. For social aspects, S1 has a 90.7% probability of ranking first. S2 ranks first in environmental sustainability with a 100% probability and ranks first in technology sustainability with a 71.1% probability.

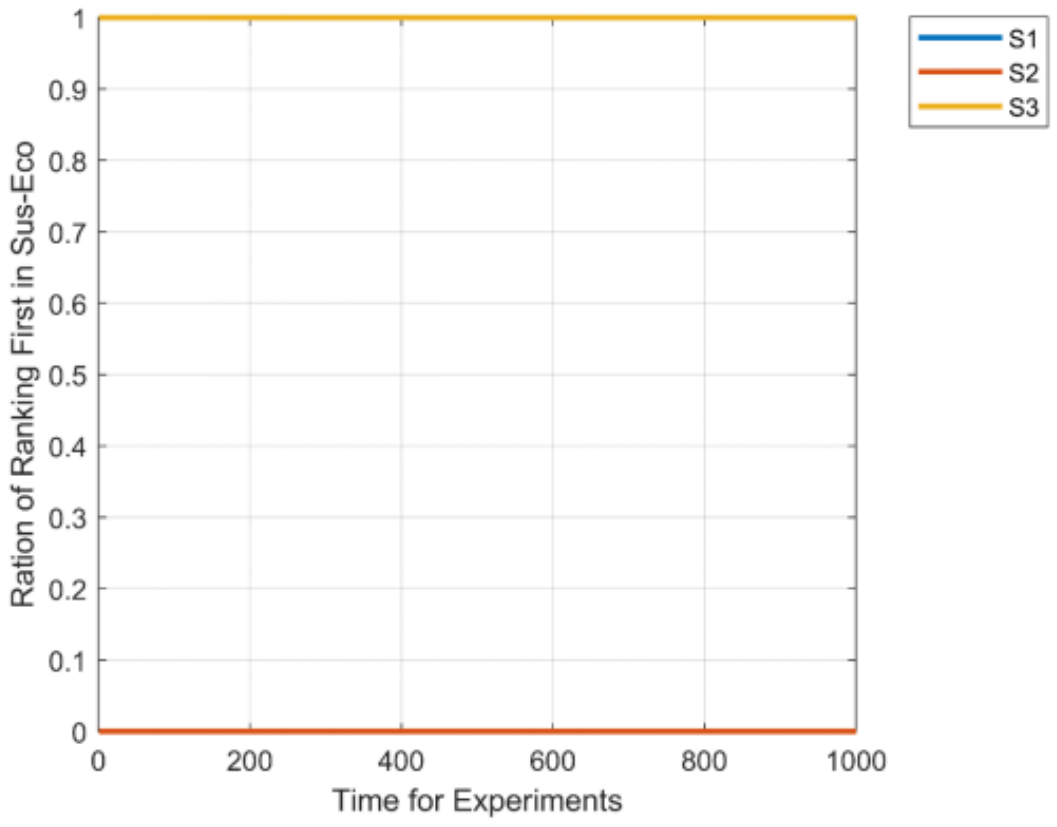

(a) Probability of ranking as first for Economic in Sustainability ( min weight 0.1)

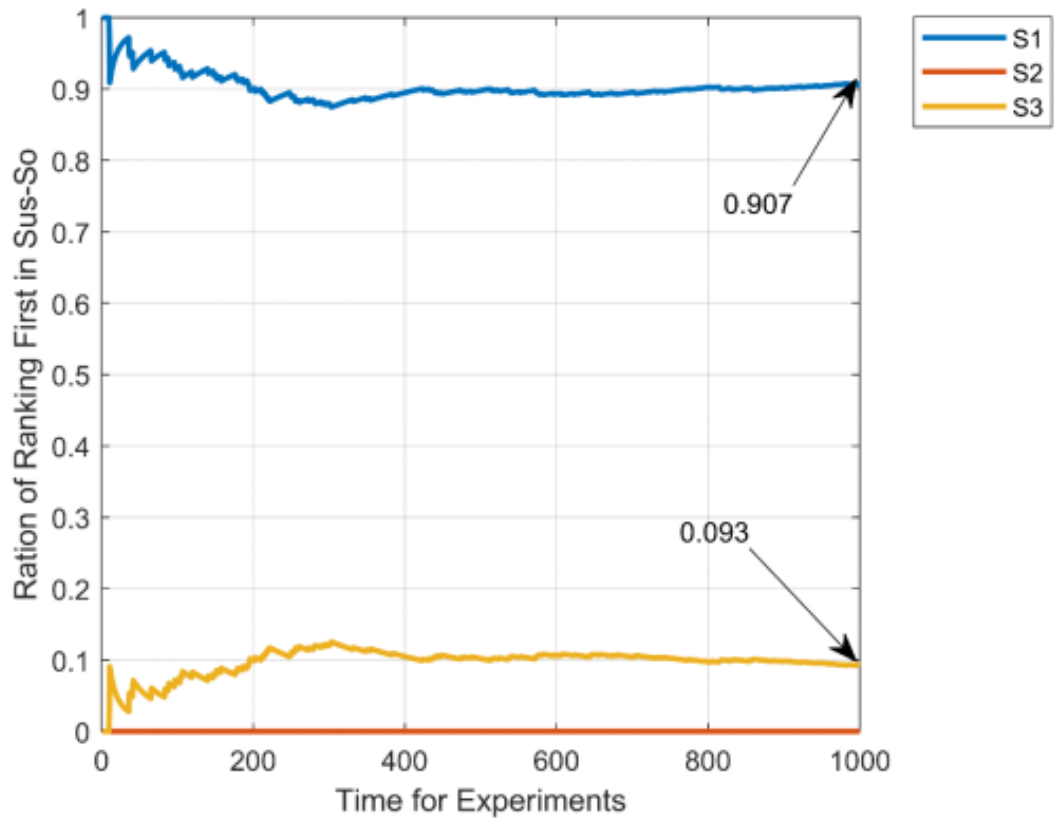

(b) Probability of ranking as first for Social in Sustainability (min weight 0.1)

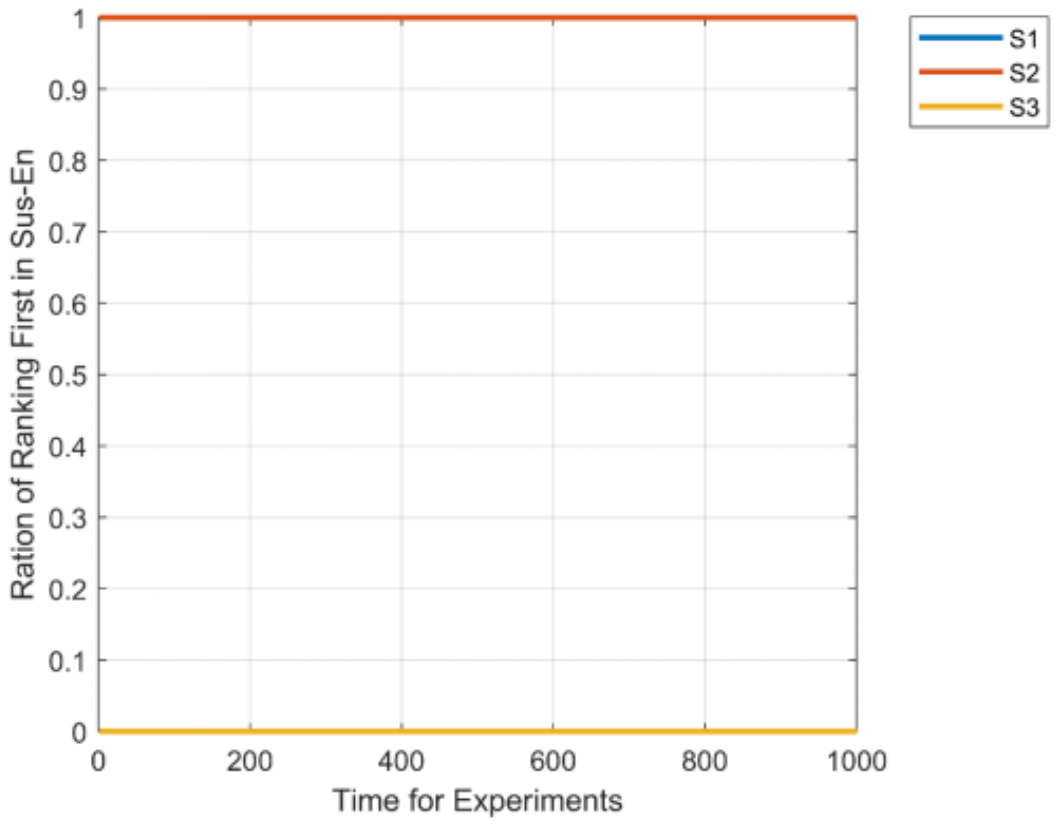

(c) Probability of ranking as first for Environment in Sustainability (min weight 0.1)

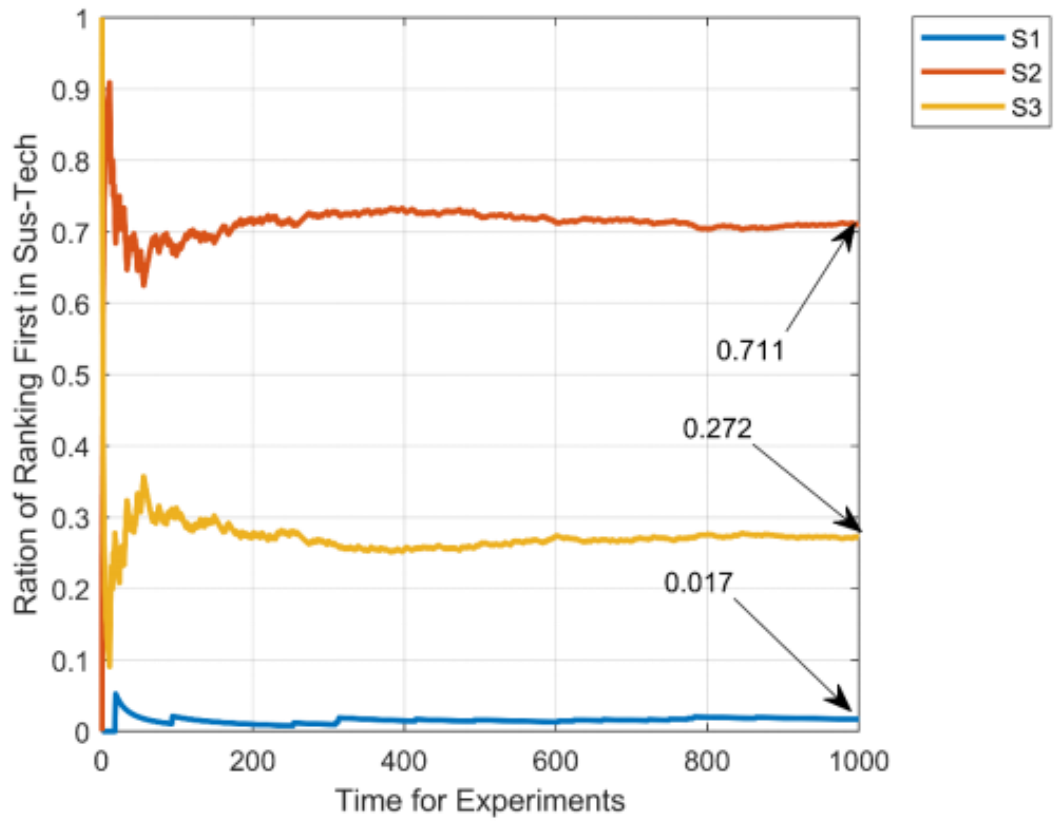

(d) Probability of ranking as first for Technology in Sustainability (Min Weight 0.1)

Figure 13 Probability of ranking as first for four requirements in Sustainability

The following figures demonstrate the more detailed information in 12 criteria. Upon analysis, it was observed that S1 obtained the highest rankings in C2, C4, C5, C7, C9, and C11 with probabilities exceeding 50%. S2, on the other hand, performed exceptionally well in C6, C8, and C10 with a probability of 100%. Finally, in C1, C3, and C12, S3 was found to excel above all other scenarios.

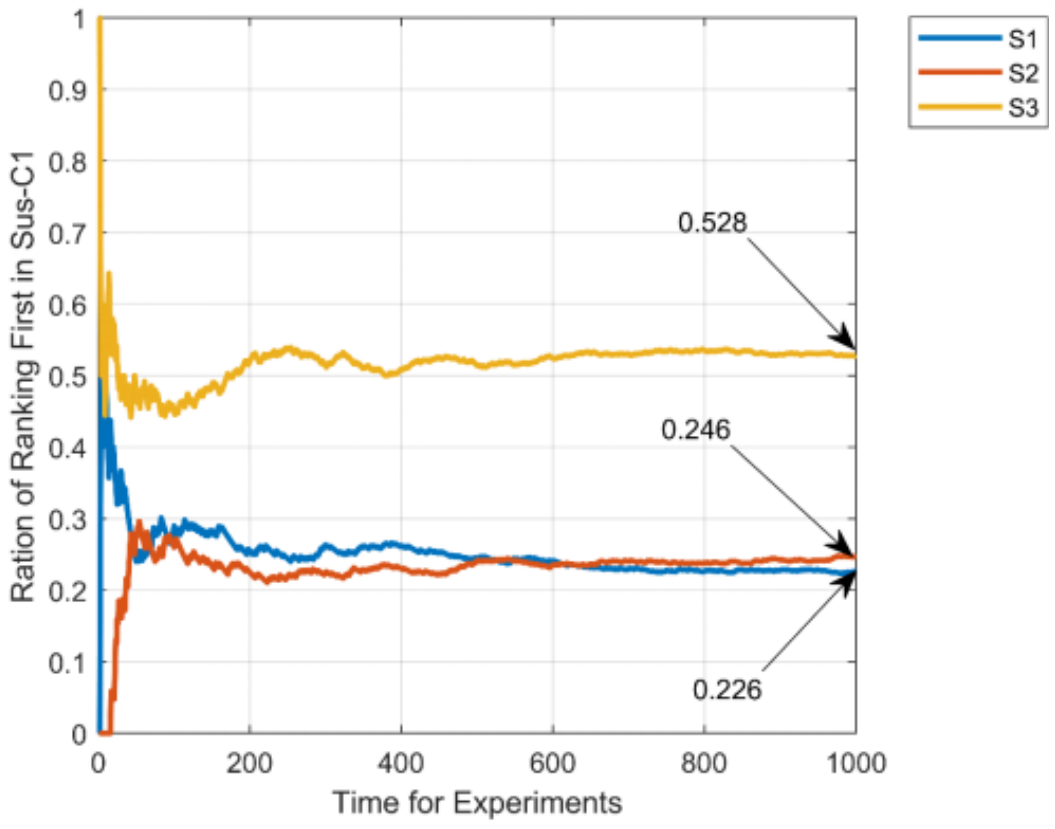

(1) Probability of ranking as first for C1 in sustainability

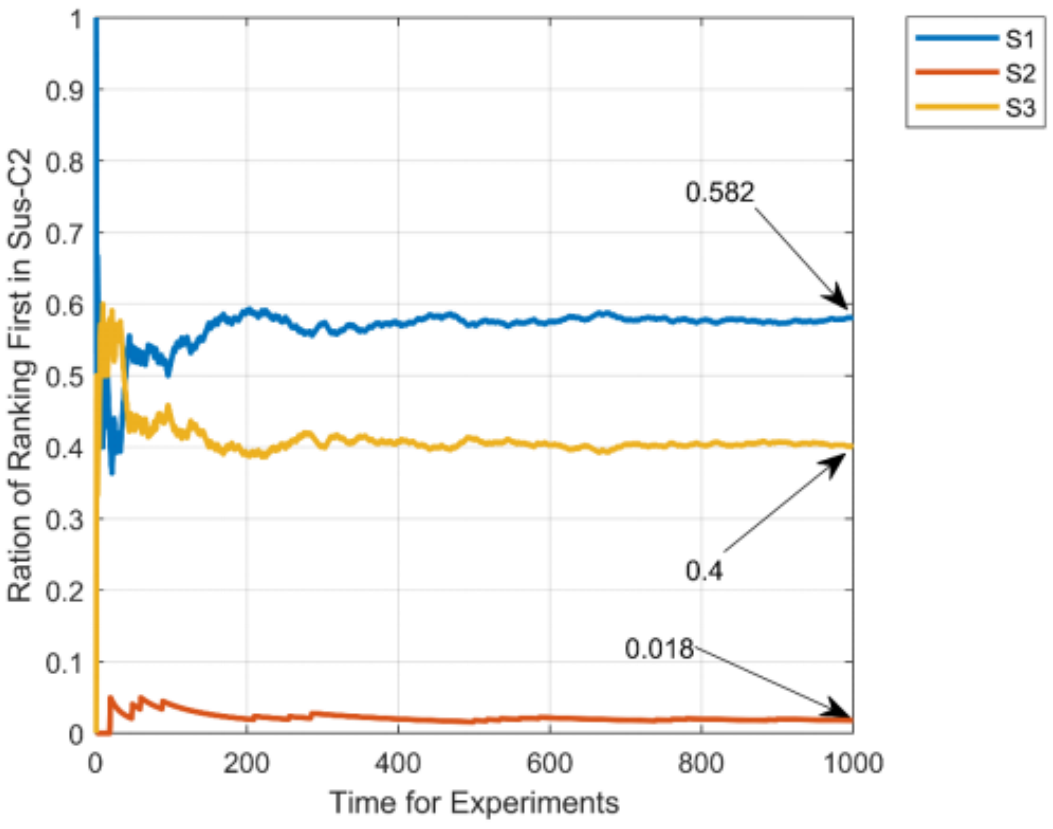

(2) Probability of ranking as first for C2 in sustainability

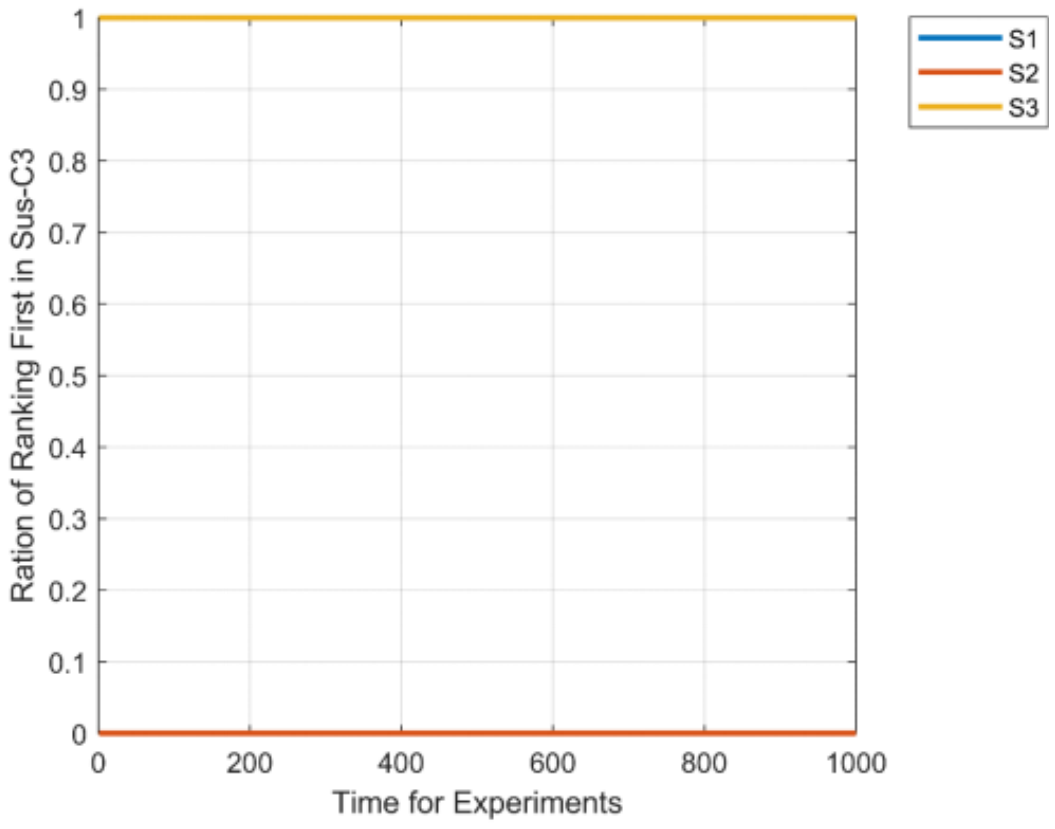

(3) Probability of ranking as first for C3 in sustainability

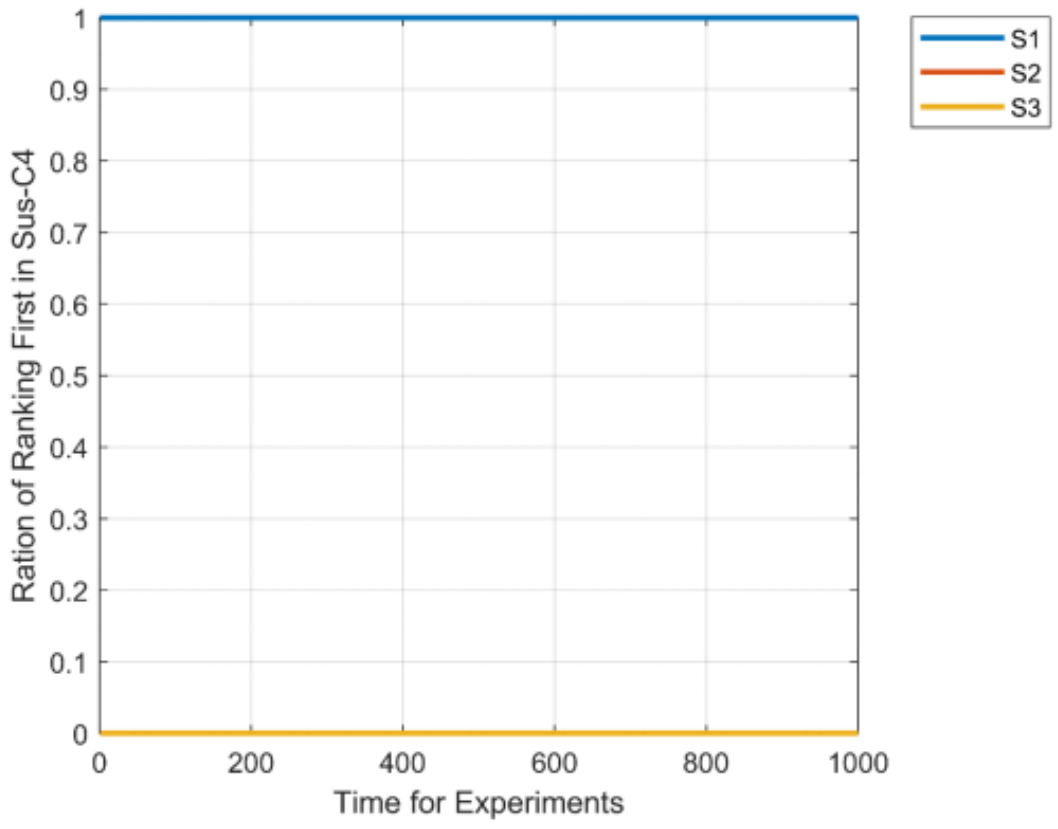

(4) Probability of ranking as first for C4 in sustainability

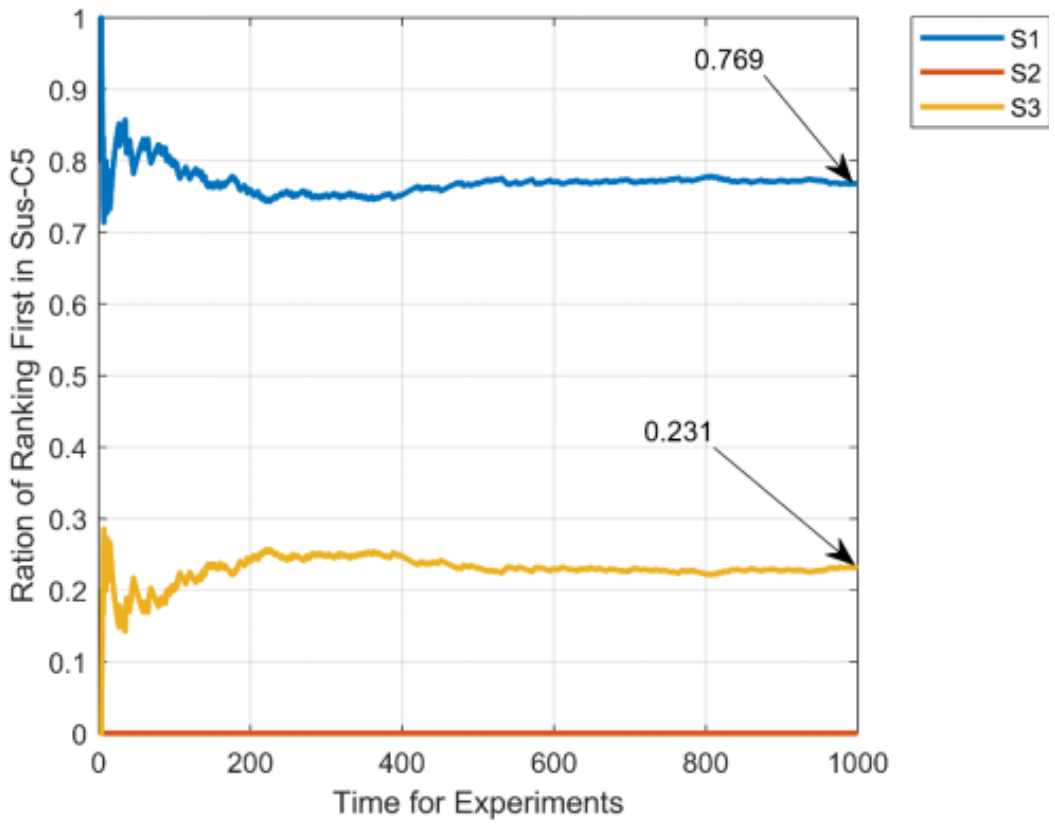

(5) Probability of ranking as first for C5 in sustainability

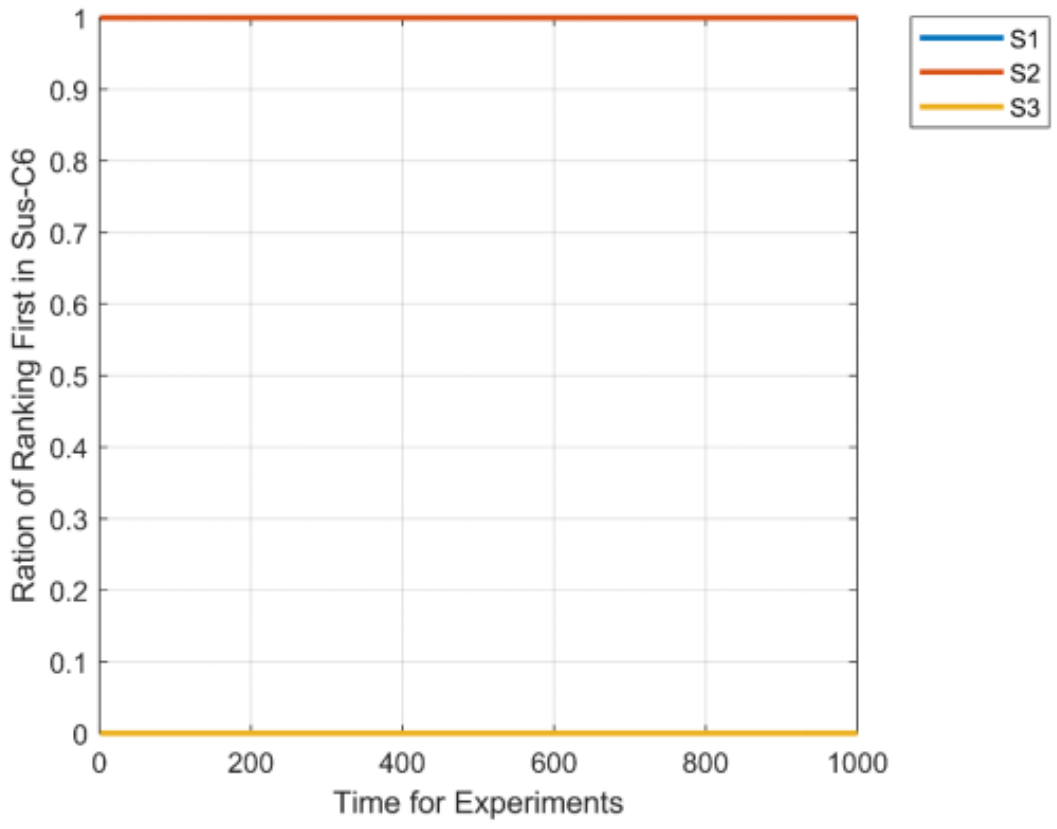

(6) Probability of ranking as first for C6 in sustainability

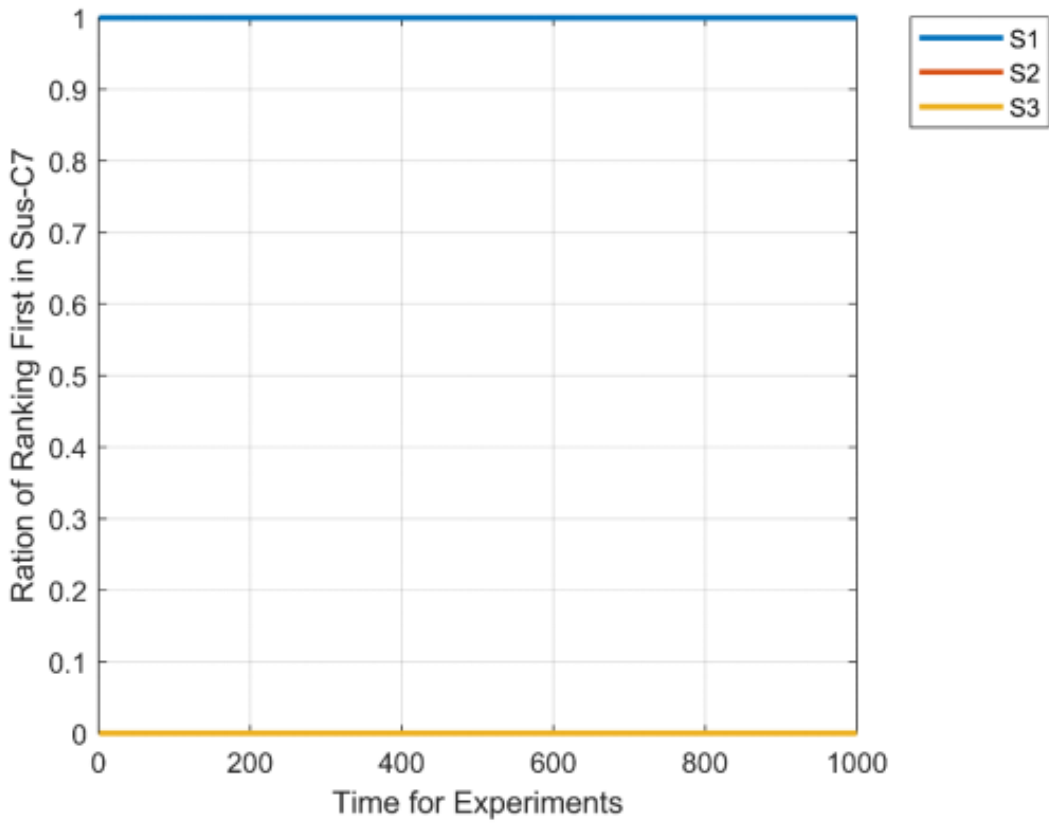

(7) Probability of ranking as first for C7 in sustainability

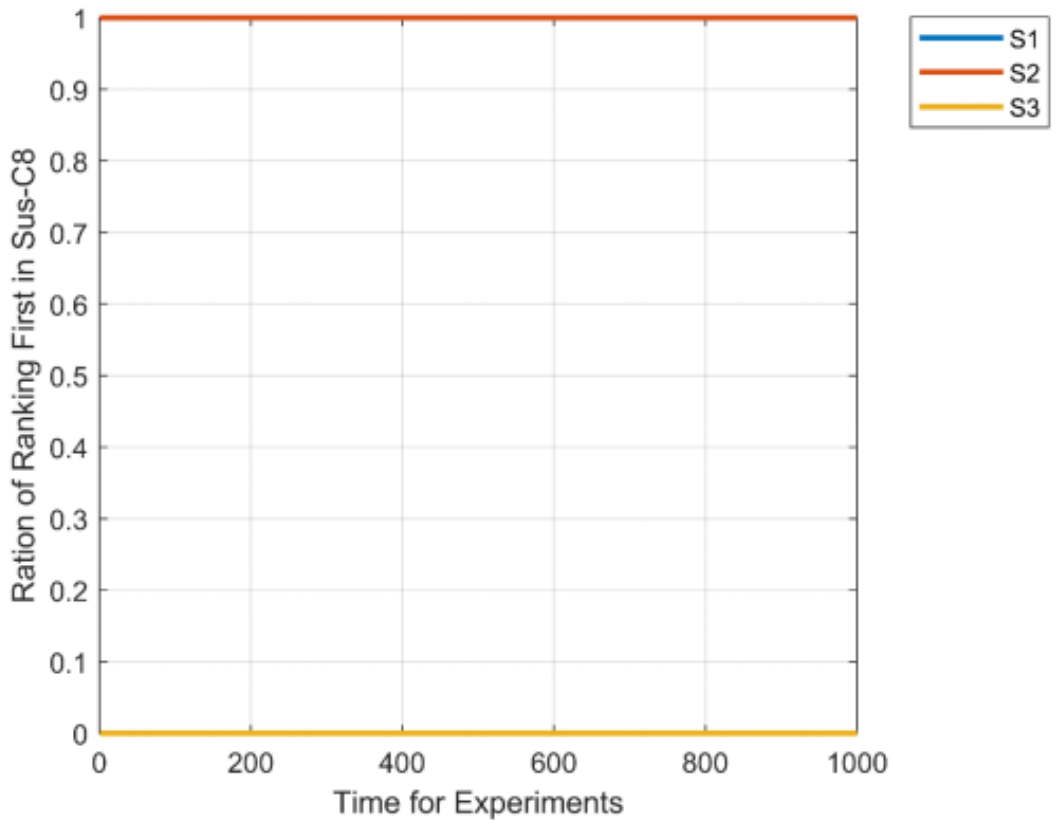

(8) Probability of ranking as first for C8 in sustainability

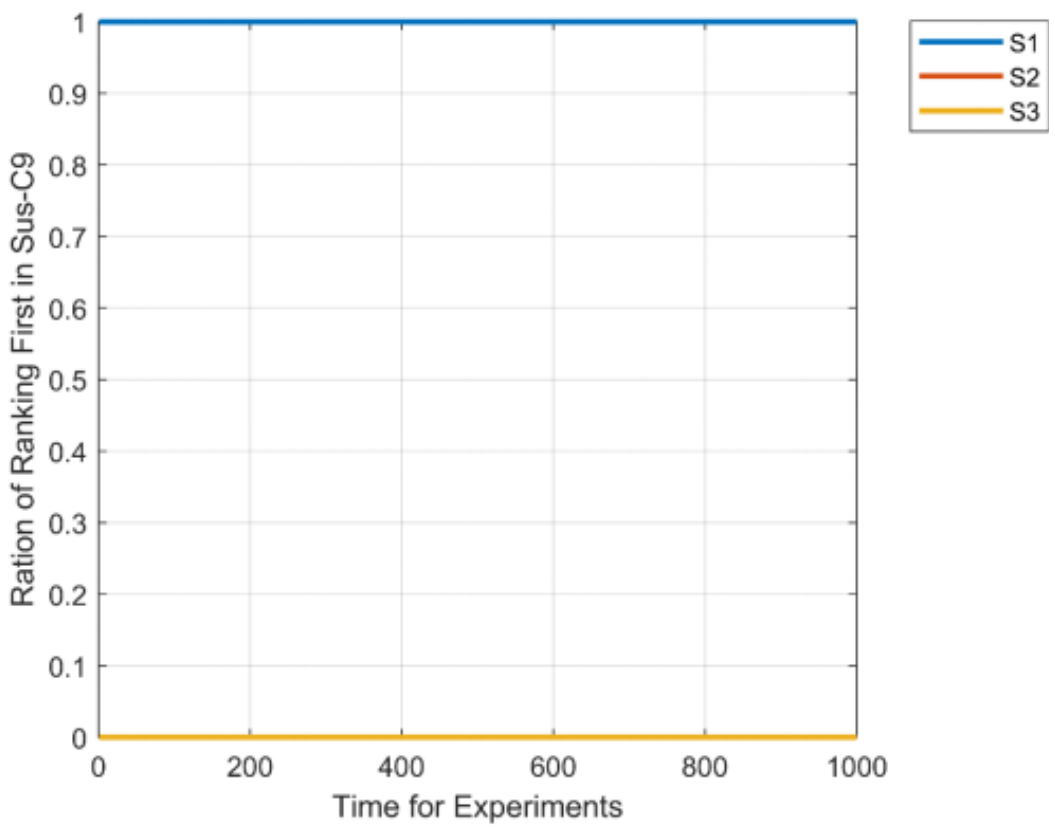

(9) Probability of ranking as first for C9 in sustainability

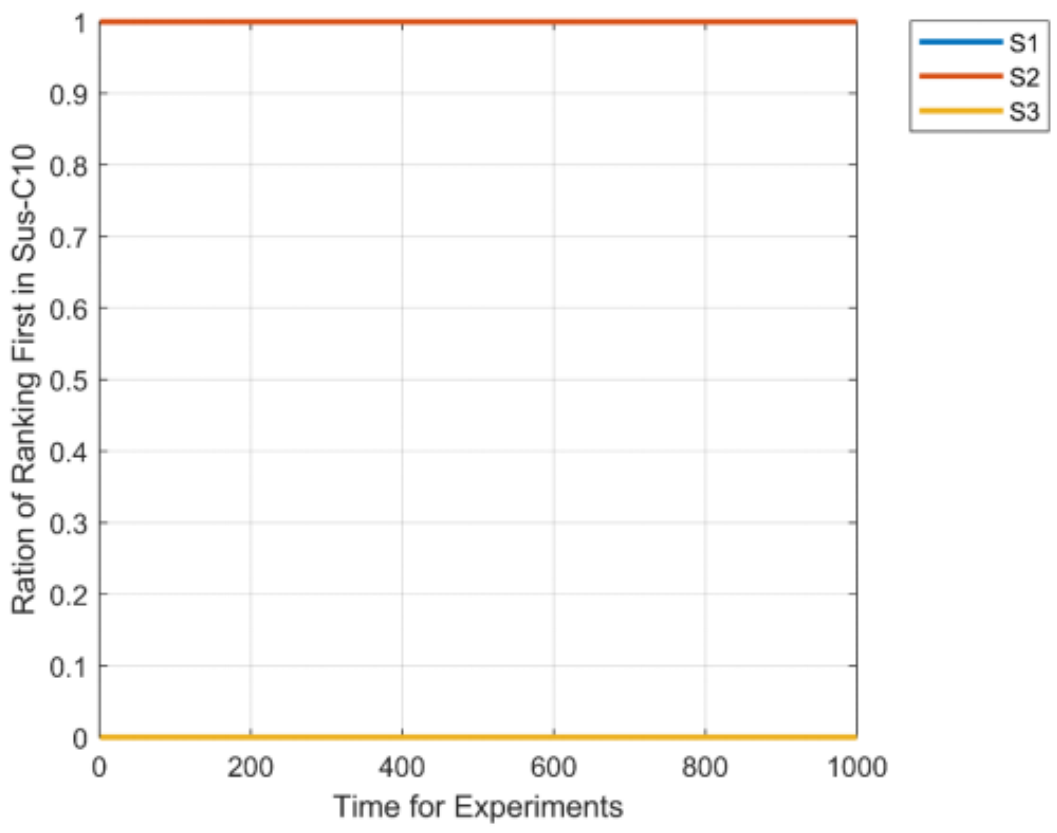

(10) Probability of ranking as first for C10 in sustainability

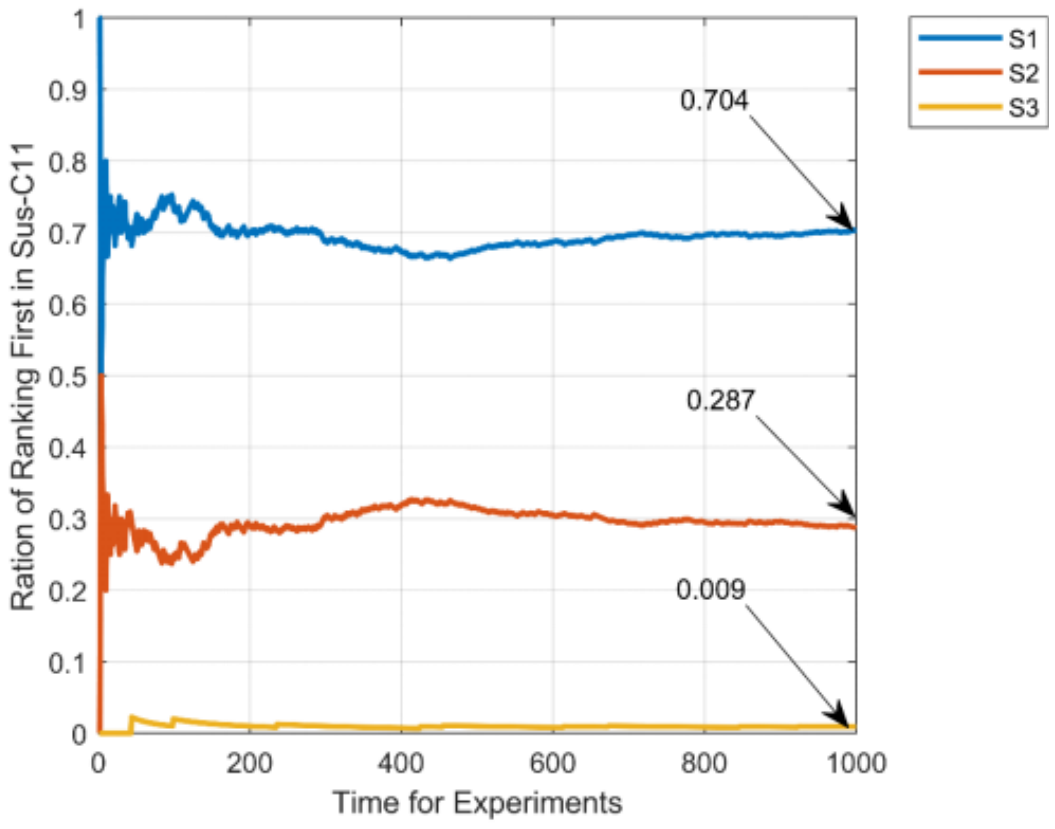

(11) Probability of ranking as first for C11 in sustainability

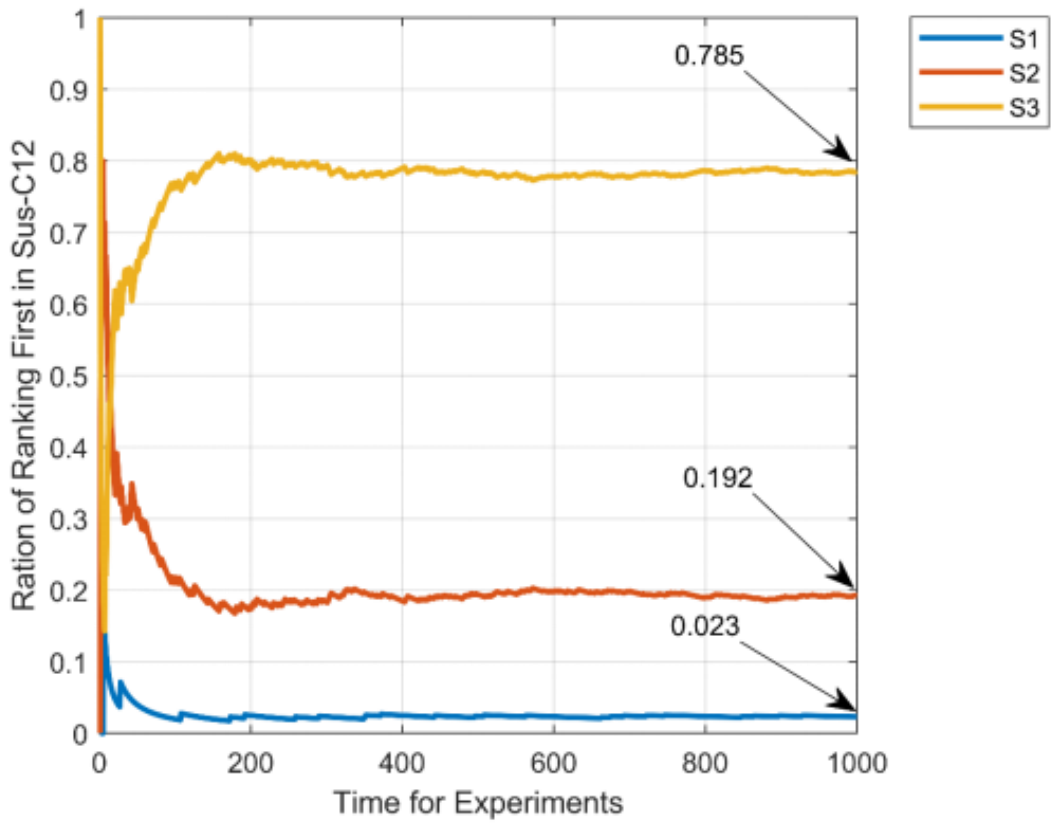

(12) Probability of ranking as first for C12 in sustainability

Figure 14 Probability of ranking as first for 12 criteria in Sustainability

## 5.3 Probability Results for Circularity

We conducted the same experiments for circularity. The Figure 15 displays the PDF and CDF representing the reliability of the simulation data. The value ranges for S1, S2, and S3 were determined to be 0.56-0.66, 0.58-0.67, and 0.61-0.68, respectively. These values are greater than those for sustainability.

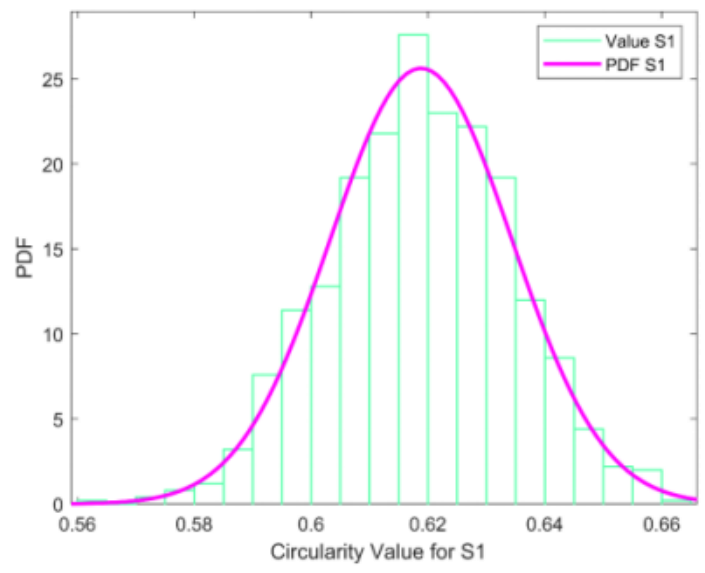

(a) PDF for the circularity value of S1

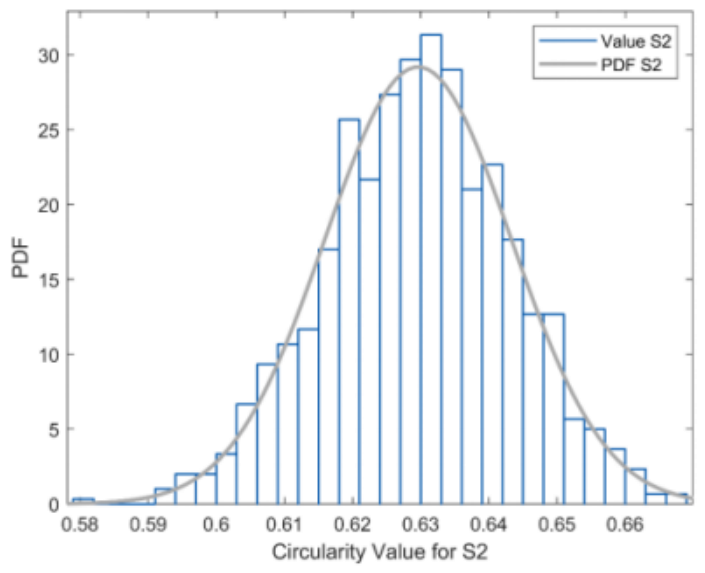

(b) PDF for the circularity value of S2

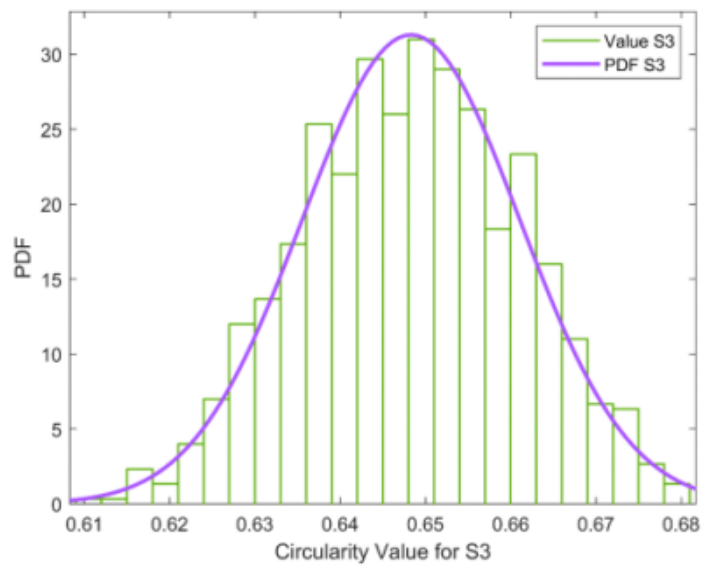

(c) PDF for the circularity value of S3

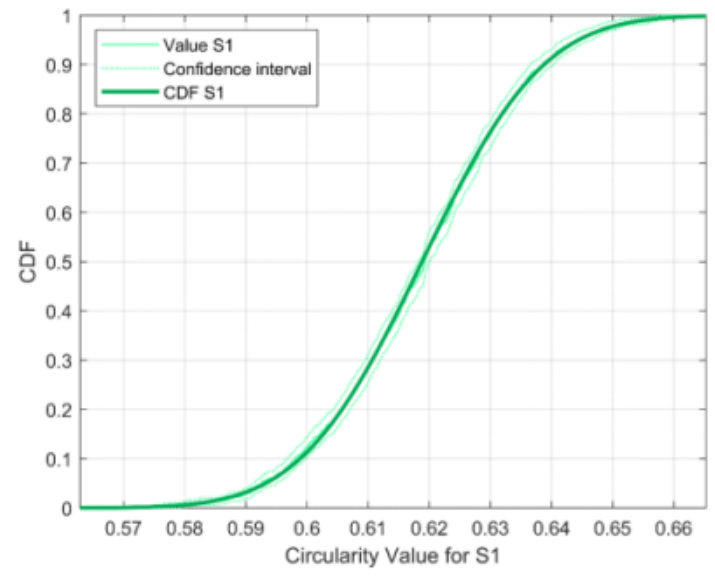

(d) CDF for the circularity value of S1

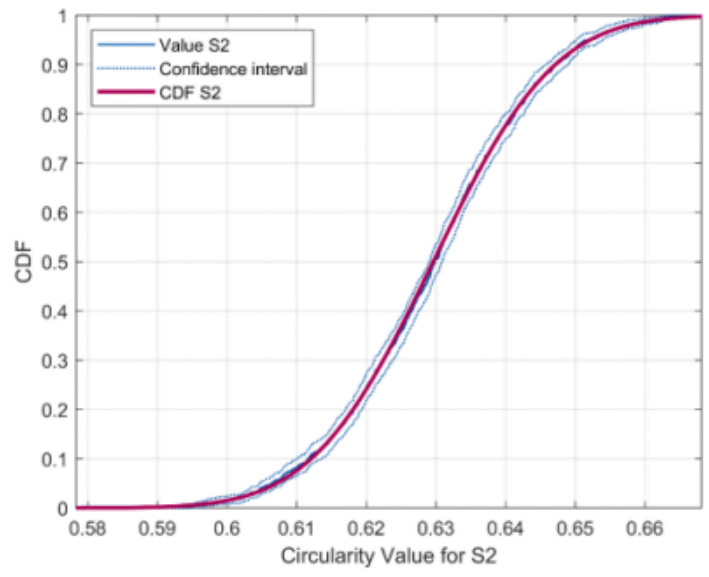

(e) CDF for the circularity value of S2

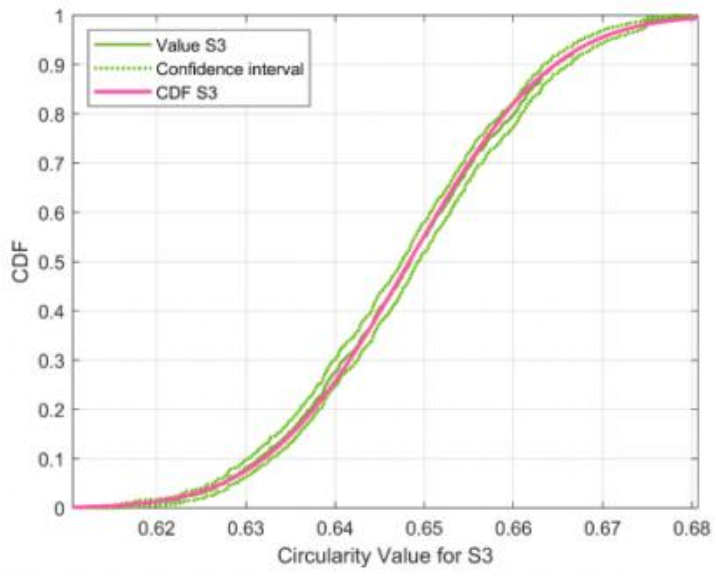

(f) CDF for the circularity value of S3

Figure 15 PDF and CDF for circularity value for three scenarios

In the Figure 16 provided, the average values of circularity levels for three different MMC product scenarios are presented. Out of these scenarios, Scenario 3 has the highest sustainability level at 0.6483, followed by Scenario 2 with 0.6296. The lowest sustainability level is attributed to Scenario 1, which stands at 0.6188.

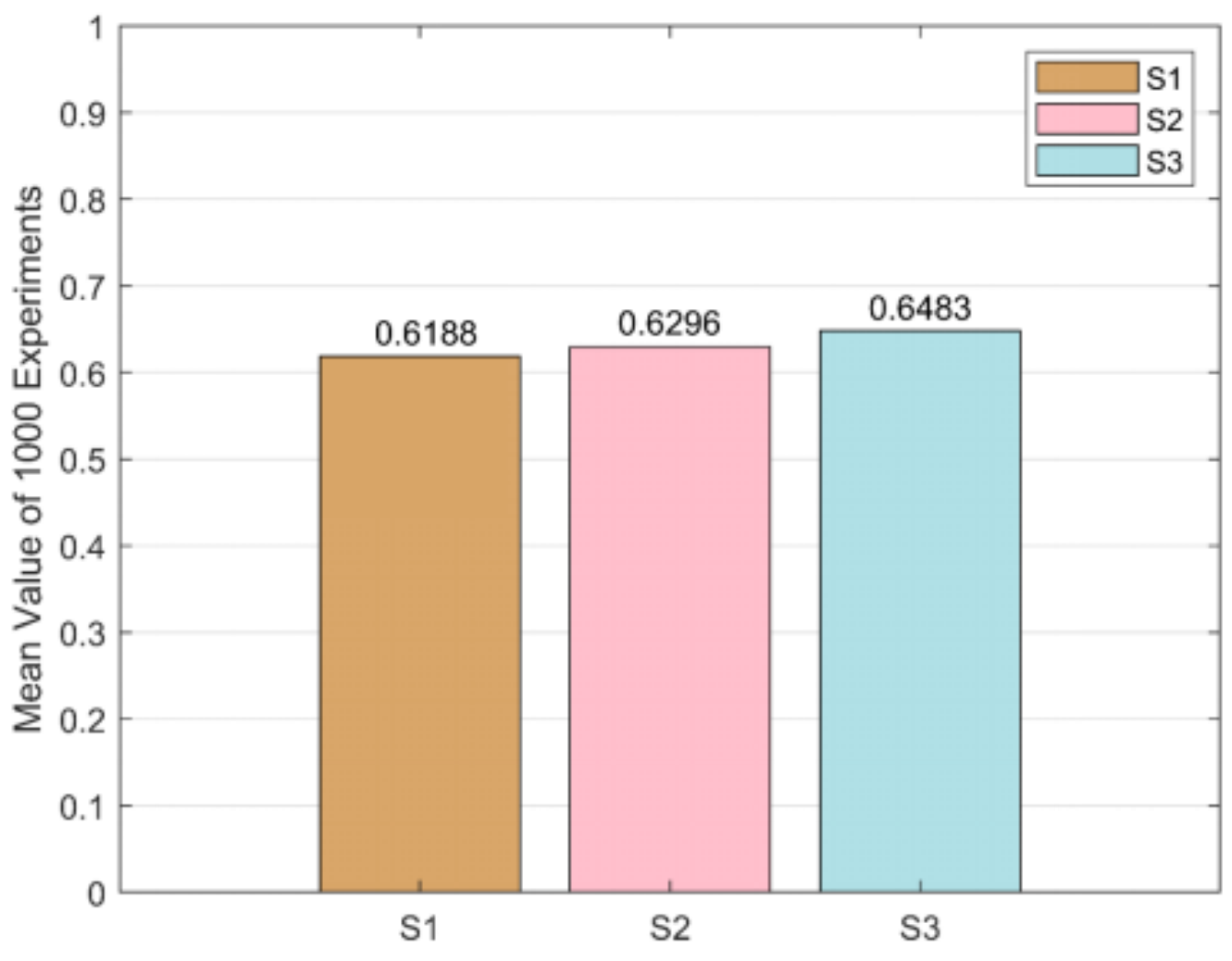

Figure 16 Mean Value of Circularity for Three Scenarios

Based on the analysis of circularity performances across four different requirement levels, it is evident that S1 stands out in terms of social and environmental performance. Meanwhile, S2 excels in technology performance, and S3 performs best in the economic aspect. The social sector is the most significant contributor to the overall sustainability value, followed by technology, environment, and economics, in descending order.

The detailed analysis of 11 criteria in the figure shows that S1 performs best in C1 for the economy, while S2 and S3 have the highest scores in C2 and C3, respectively. For the social aspect, S1 scores the highest in C4 criteria, and S3 performs similarly in C5. Regarding the environment, S2 and S3 have the same score for C6, while S1 has the best result in C7. All three scenarios have similar performances in C8. In terms of technology, S2 has the highest value for C9, S1 for C10, and S3 for C11.

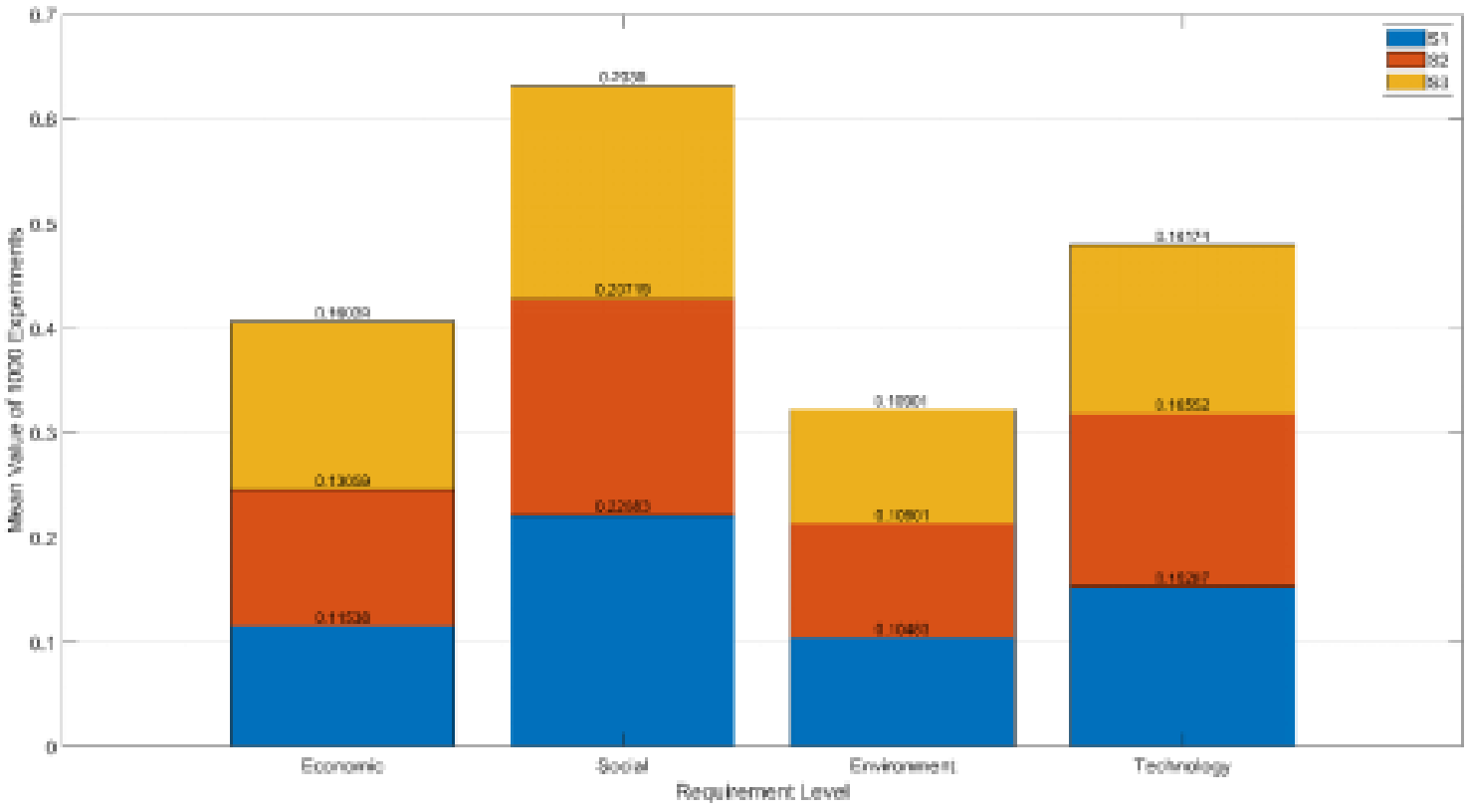

Figure 17 Stack bar for the mean value of circularity of the three scenarios in four requirements

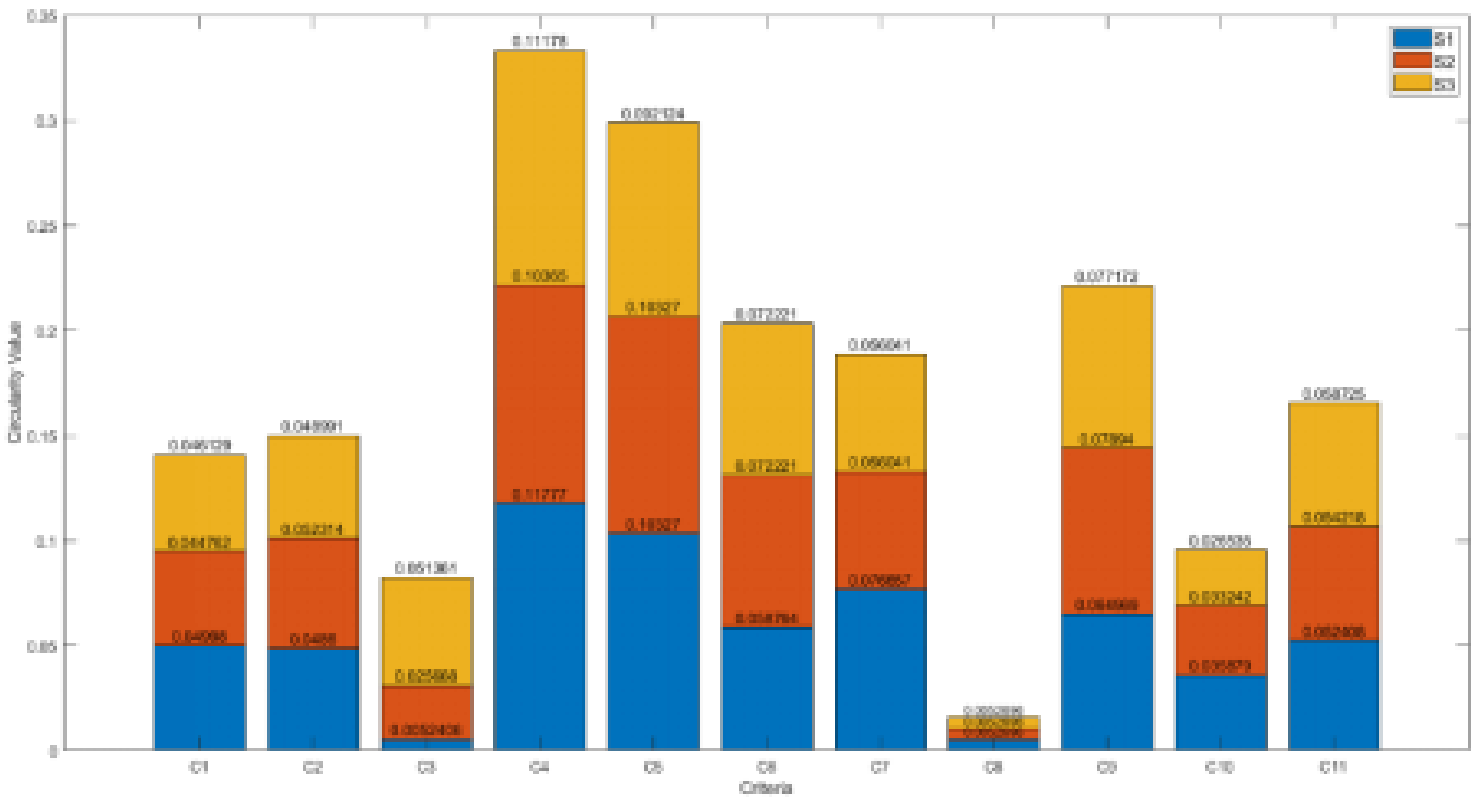

Figure 18 Stack bar for the mean value of circularity comparison between 12 criteria

After analyzing the performance of 4 different requirement levels and 12 criteria levels for three scenarios displayed in the figure, it was concluded that the circularity performances for each scenario were the same as sustainability.

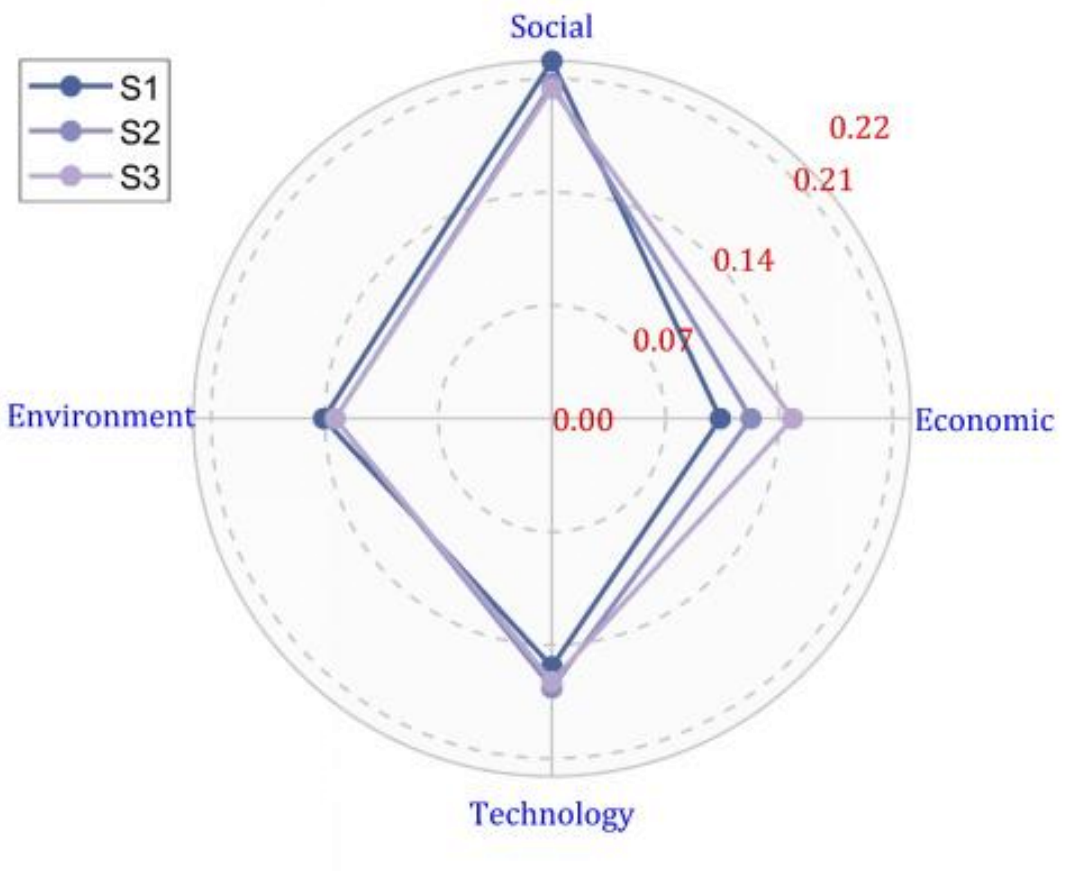

Figure 19 Radar chart for mean circularity value of four requirements of three Scenarios

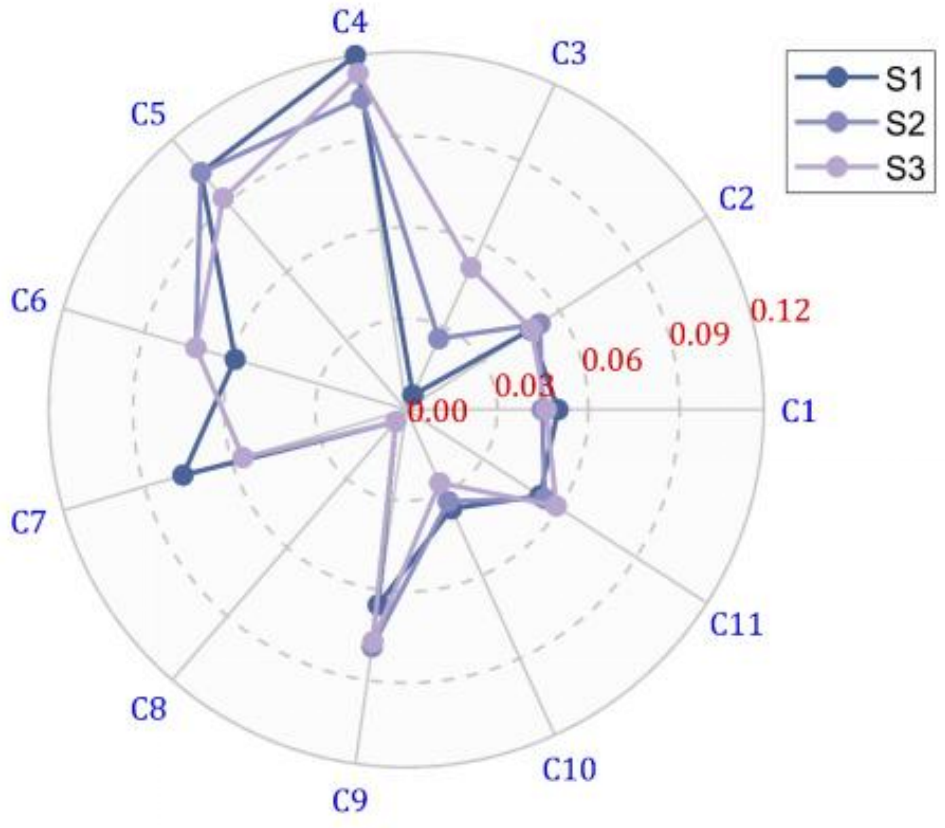

Figure 20 Radar chart for mean circularity value of 12 criteria of of three Scenarios

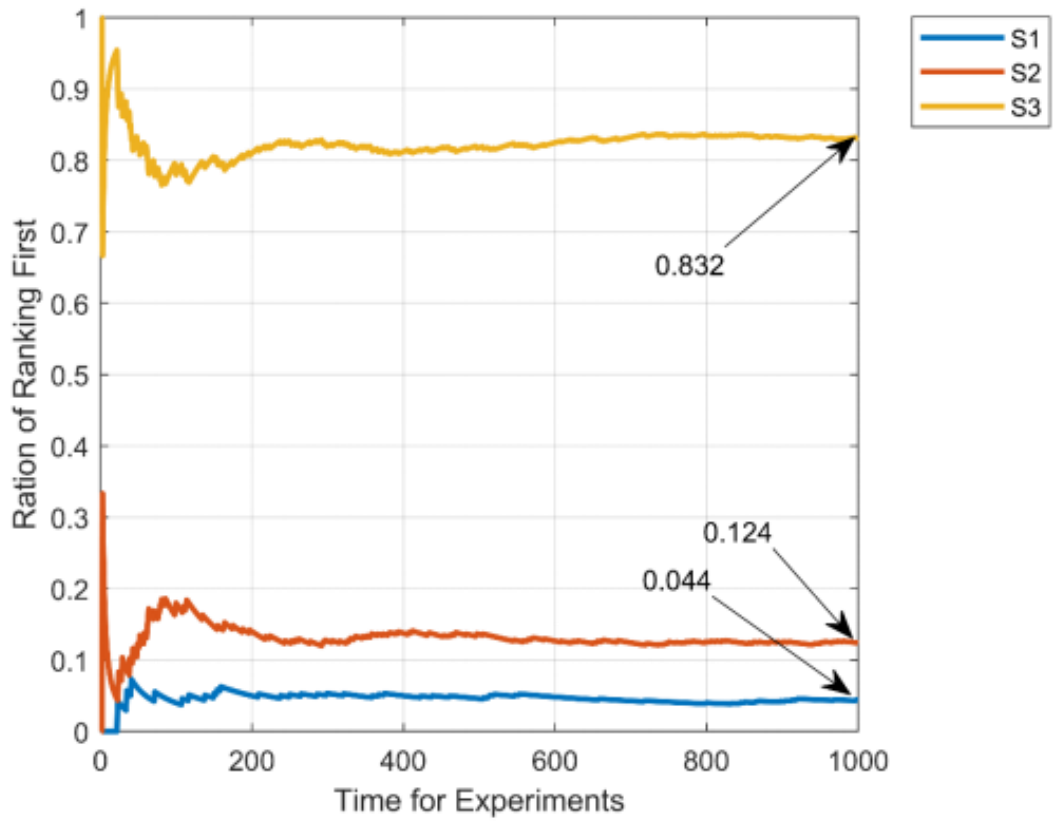

Figure 21 Probability of ranking as first for three scenarios in Circularity

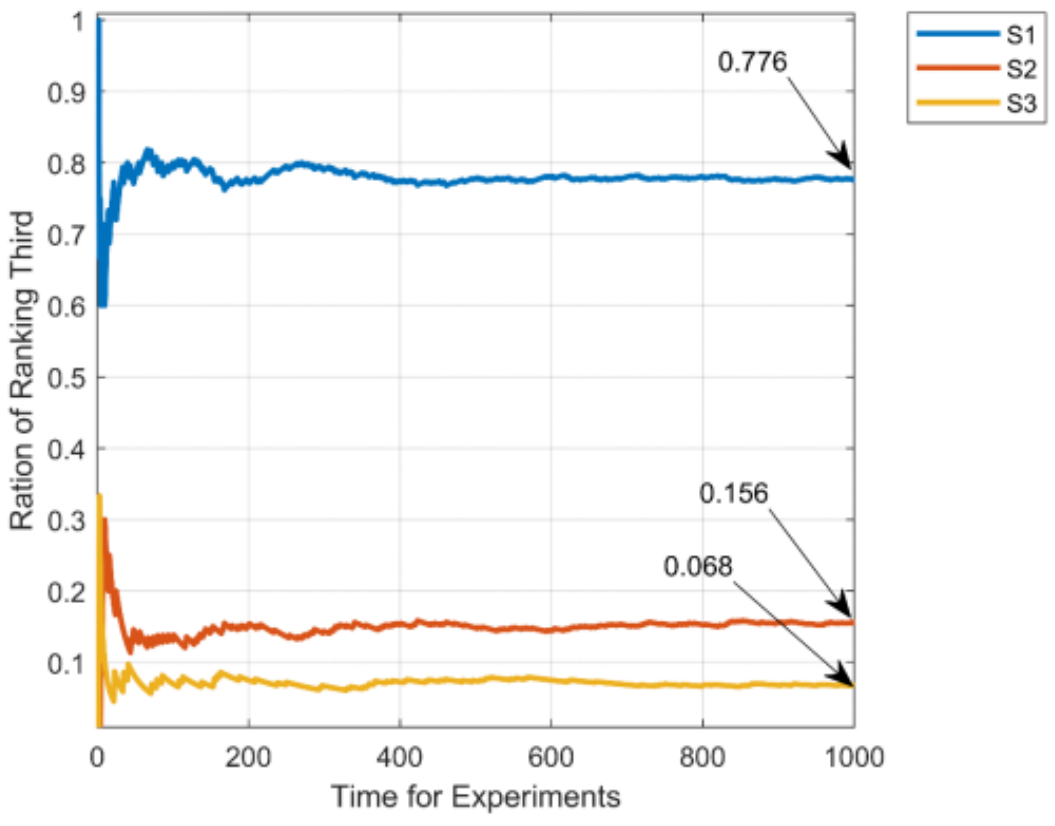

Figure 22 Probability of ranking as third for three scenarios in Circularity

The following statistics illustrate the probability of attaining the top ranking in four circularity categories. Based on the data, S3 has a 98.8% likelihood of securing the first position in economic circularity. In terms of social circualarity, S1 has a 91.7% probability of obtaining the top spot. Concerning environmental circularity, S3 is certain to come in first with a 100% probability. Finally, S2 is projected to achieve the highest ranking in technology sustainability with a 72.3% probability.

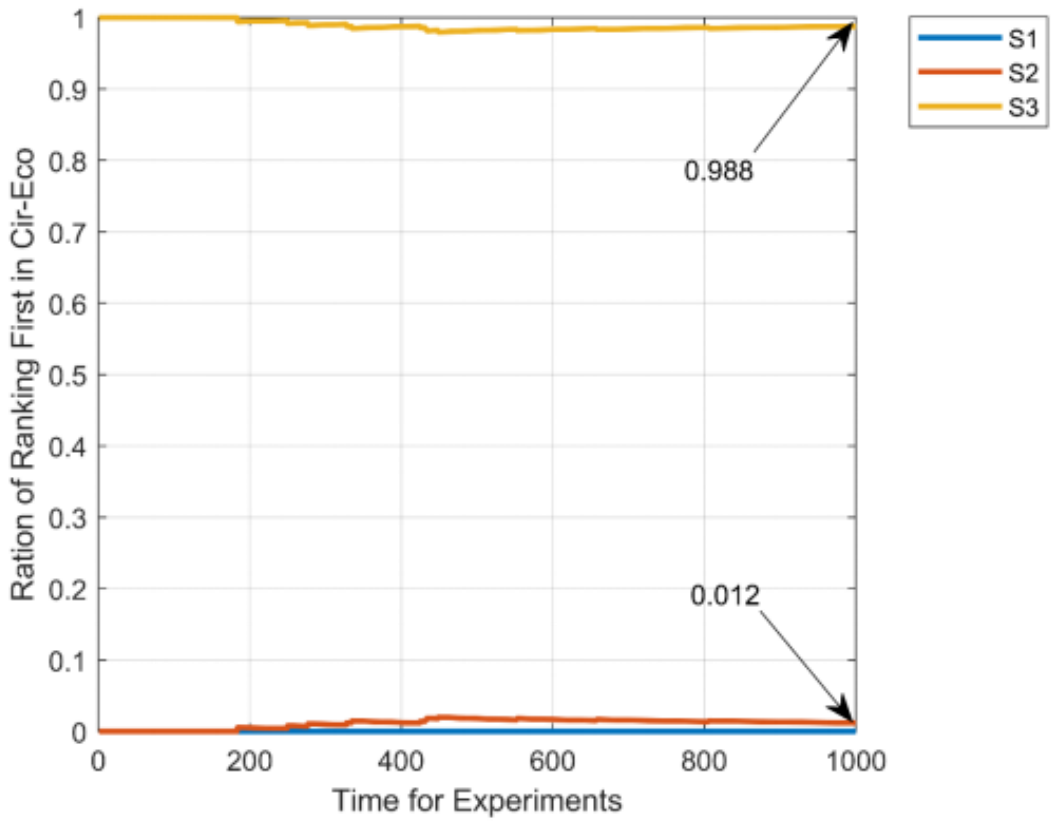

(a) Probability of ranking as first for Economic in Circularity ( min weight 0.1)

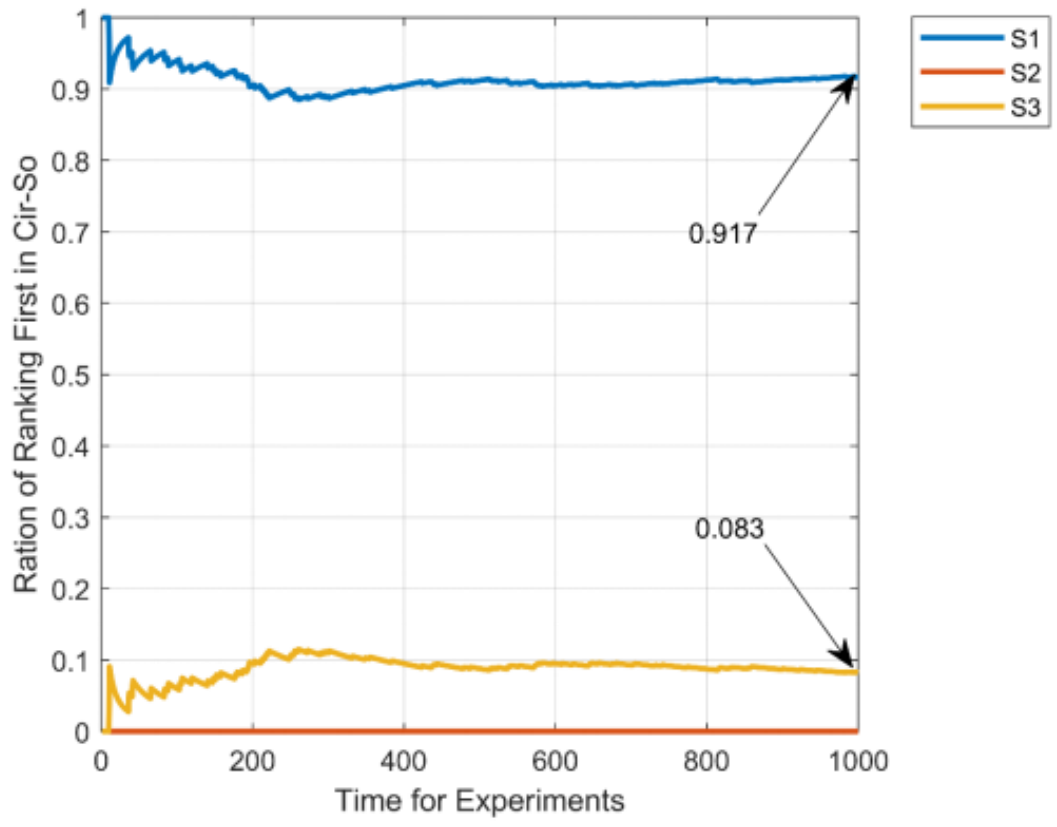

(b) Probability of ranking as first for Social in Circularity (min weight 0.1)

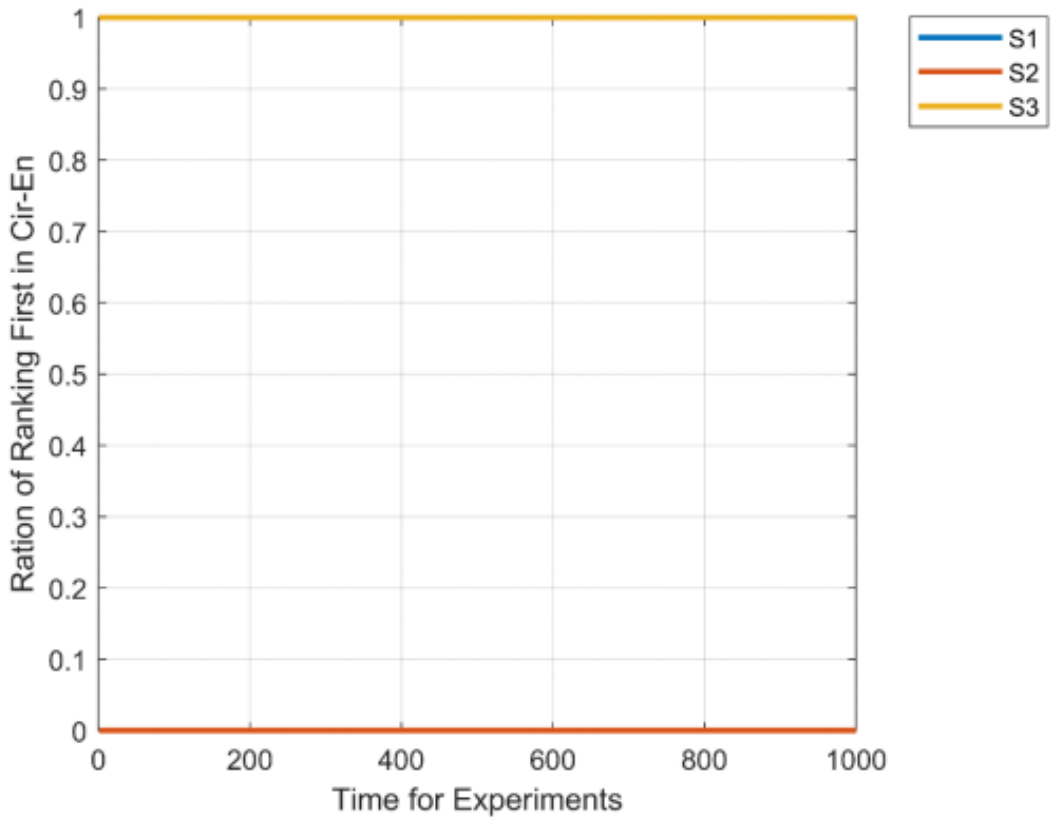

(c) Probability of ranking as first for Environment in Circularity (min weight 0.1)

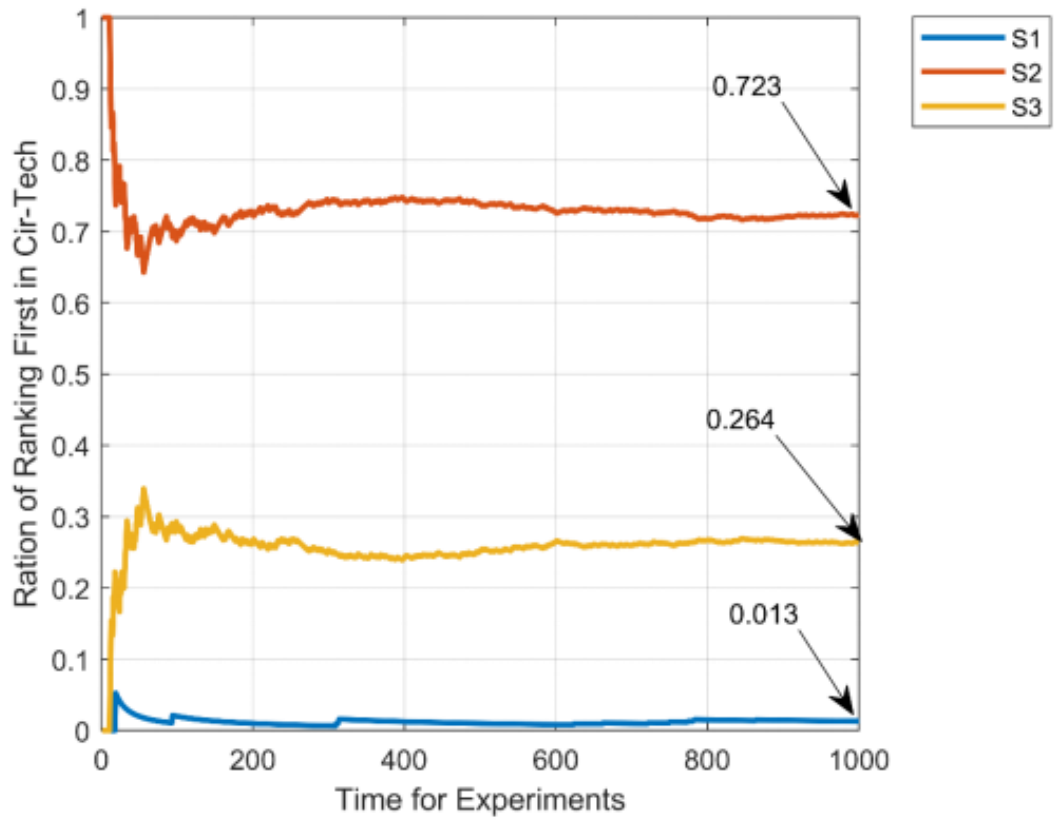

(d) Probability of ranking as first for Technology in Circularity (Min Weight 0.1)

Figure 23 Probability of ranking as first for four requirements in Circularity

The following figures demonstrate the more detailed information in 11 criteria. Upon analysis, it was observed that S1 obtained the highest rankings in C1, C4, C5, C7, C8 and C10 with probabilities exceeding 50%. S2, on the other hand, performed exceptionally well in C6 and C9 with a probability of 100%. Finally, in C2, C3, and C11, S3 was found to excel above all other scenarios.

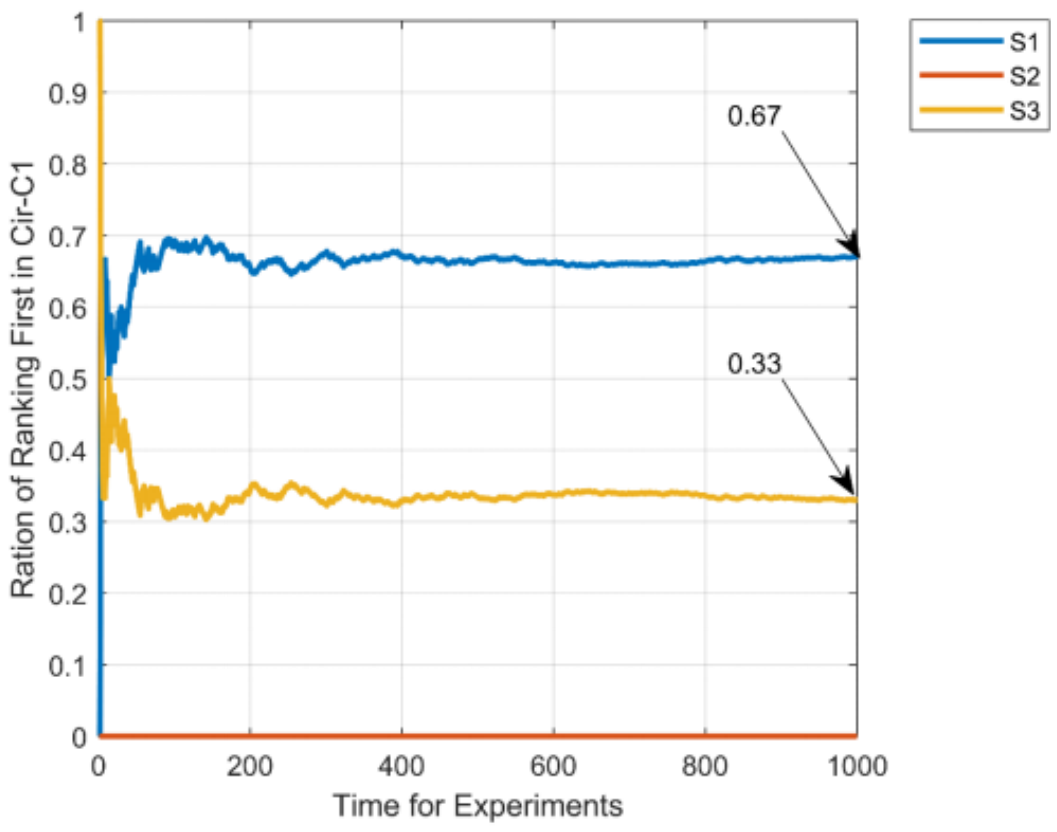

(1) Probability of ranking as first for C1 in

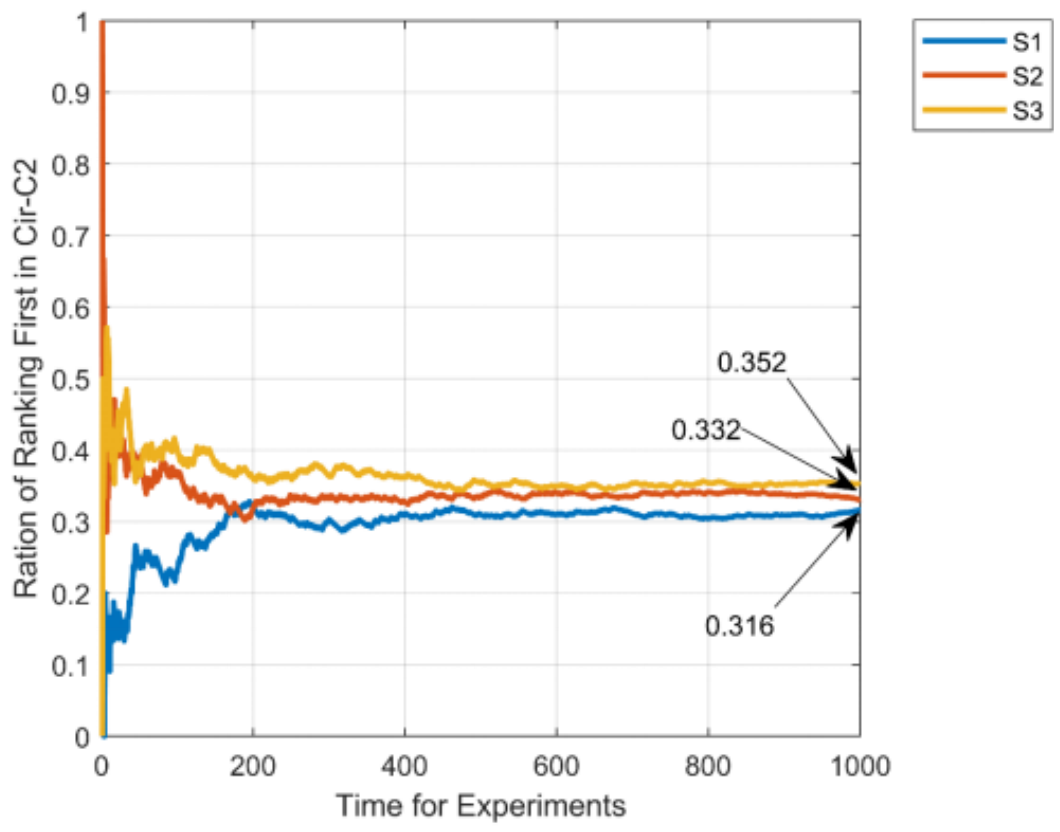

(2) Probability of ranking as first for C2 in

circularity

circularity

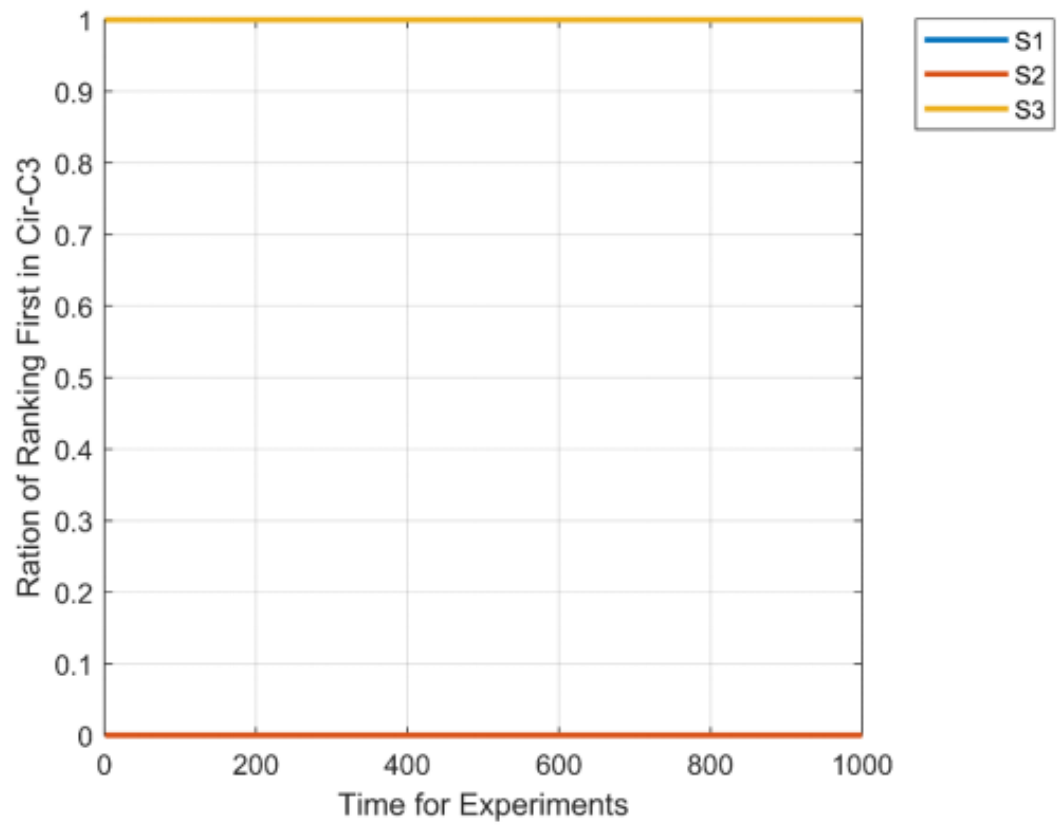

(3) Probability of ranking as first for C3 in circularity

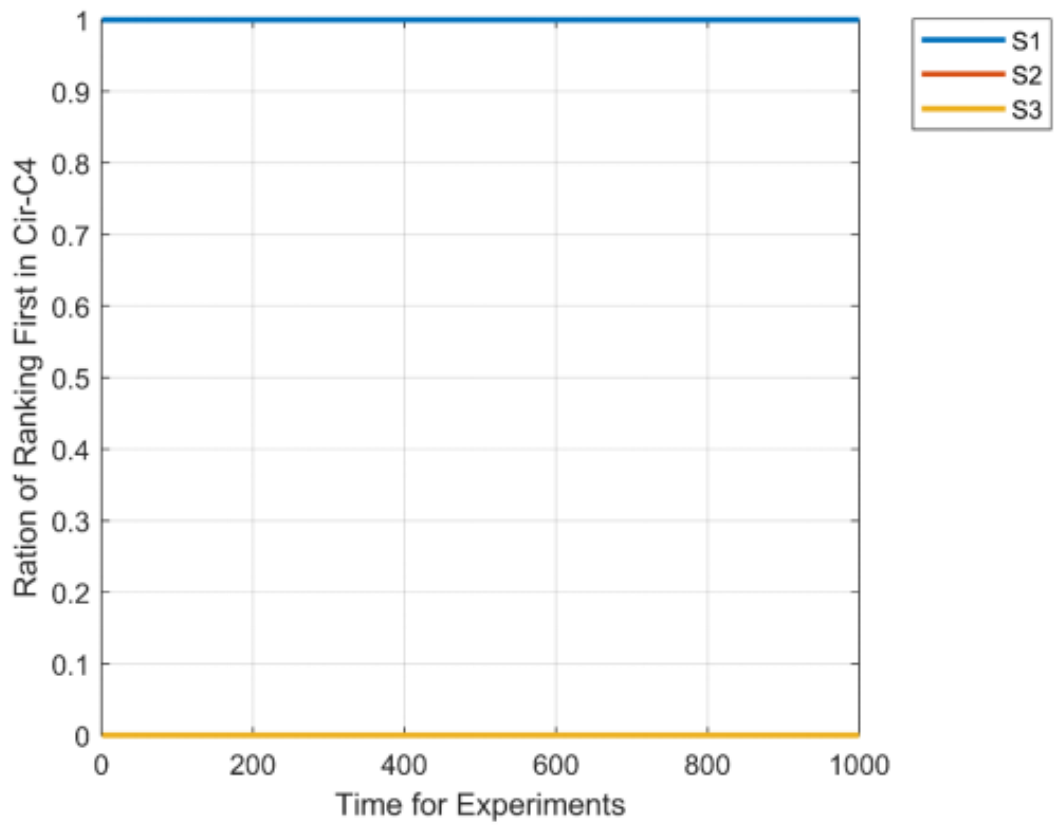

(4) Probability of ranking as first for C4 in circularity

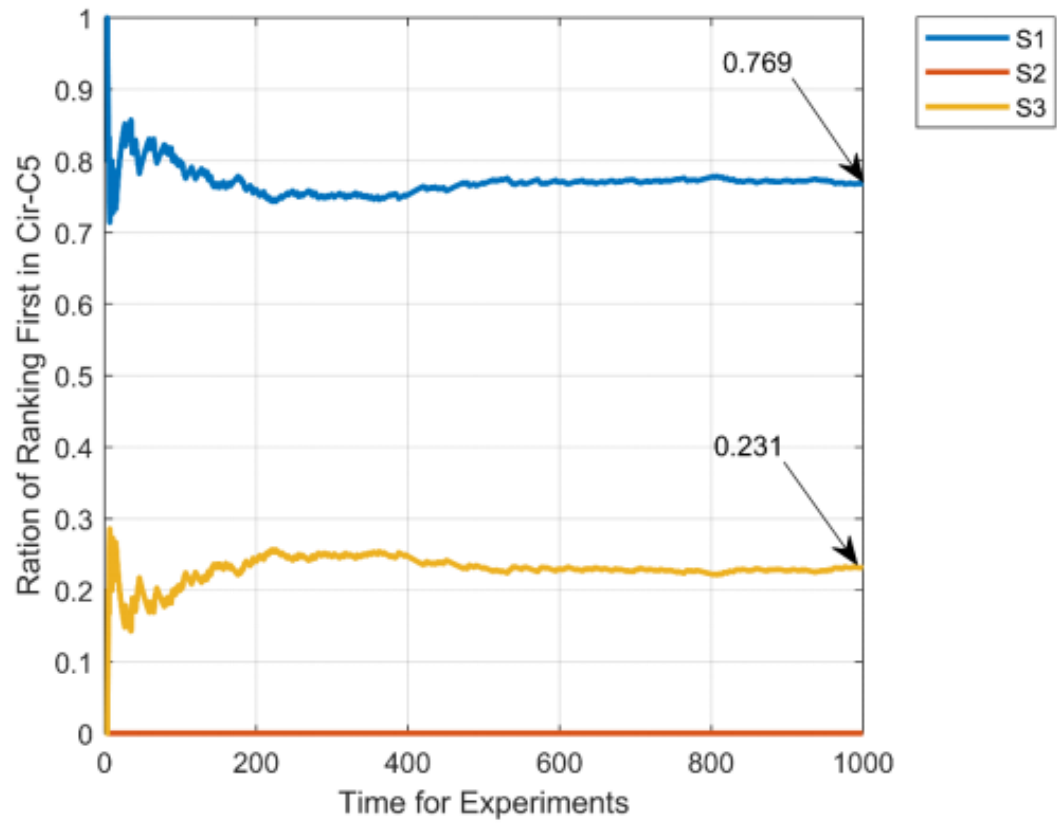

(5) Probability of ranking as first for C5 in circularity

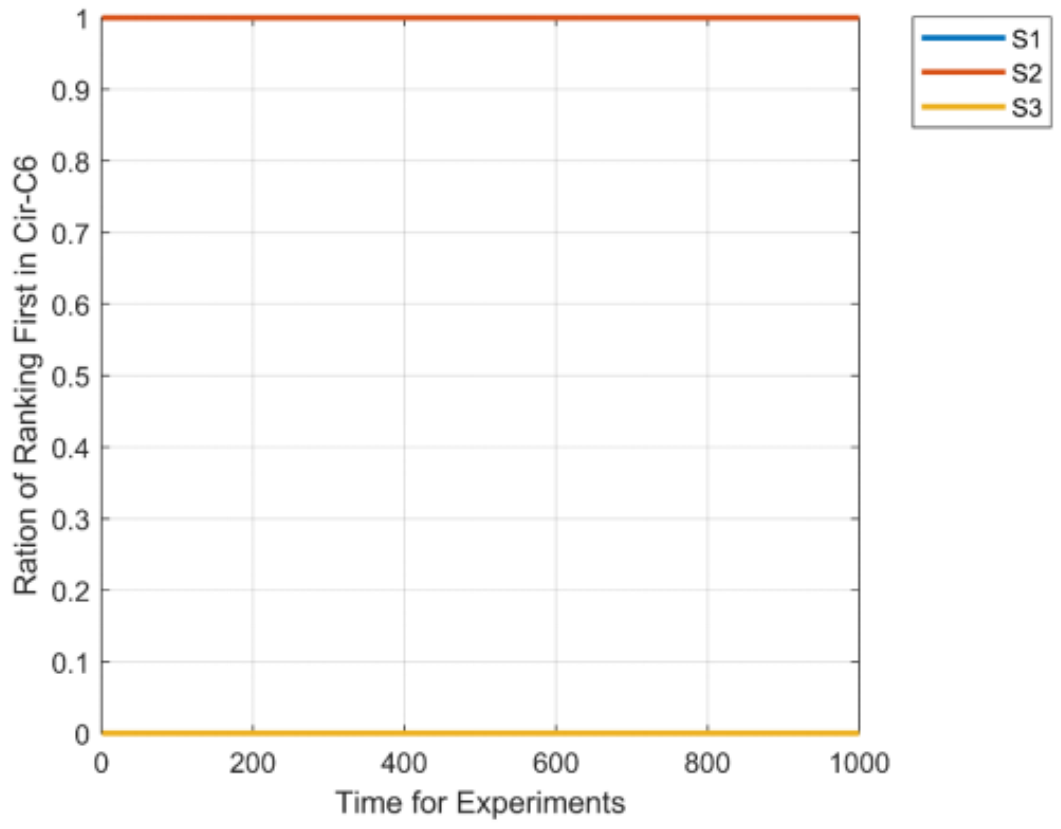

(6) Probability of ranking as first for C6 in circularity

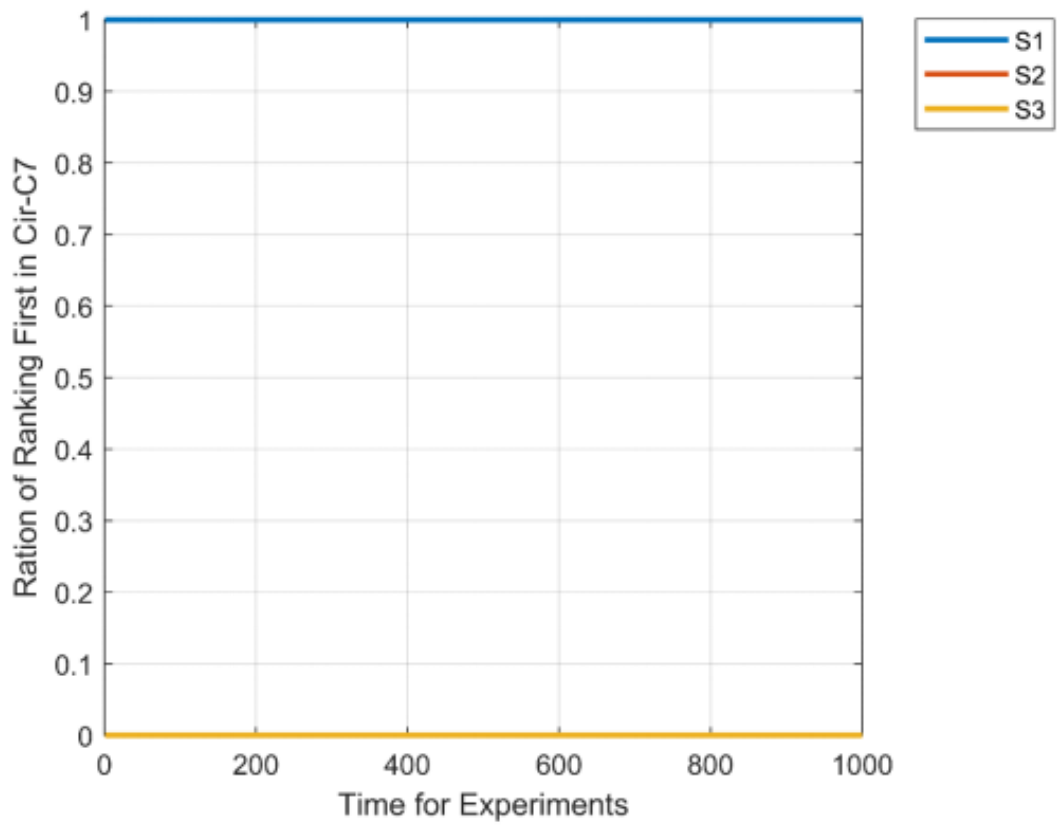

(7) Probability of ranking as first for C7 in circularity

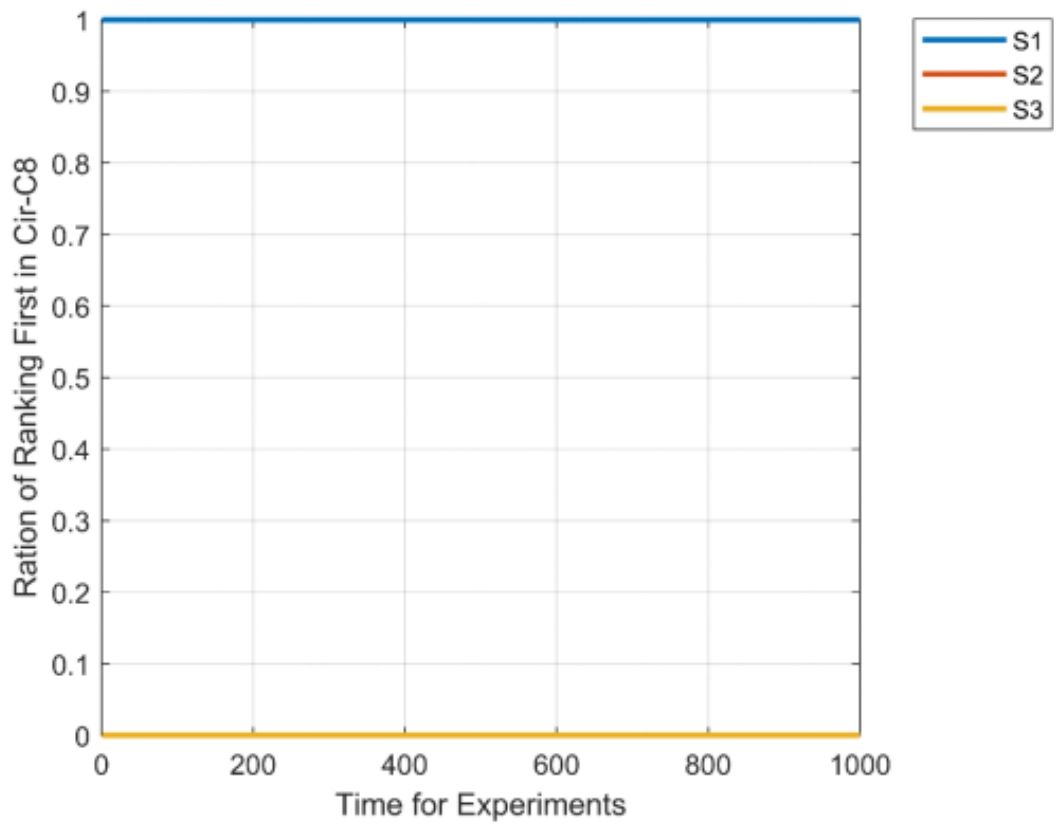

(8) Probability of ranking as first for C8 in circularity

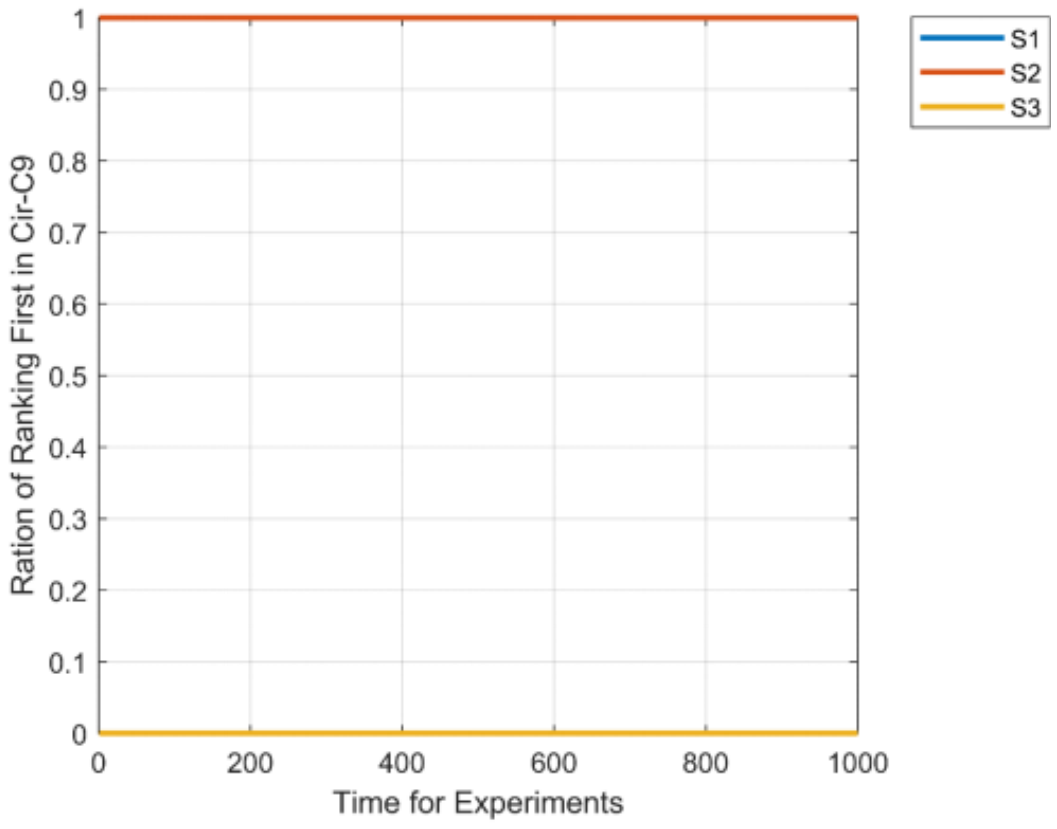

(9) Probability of ranking as first for C9 in circularity

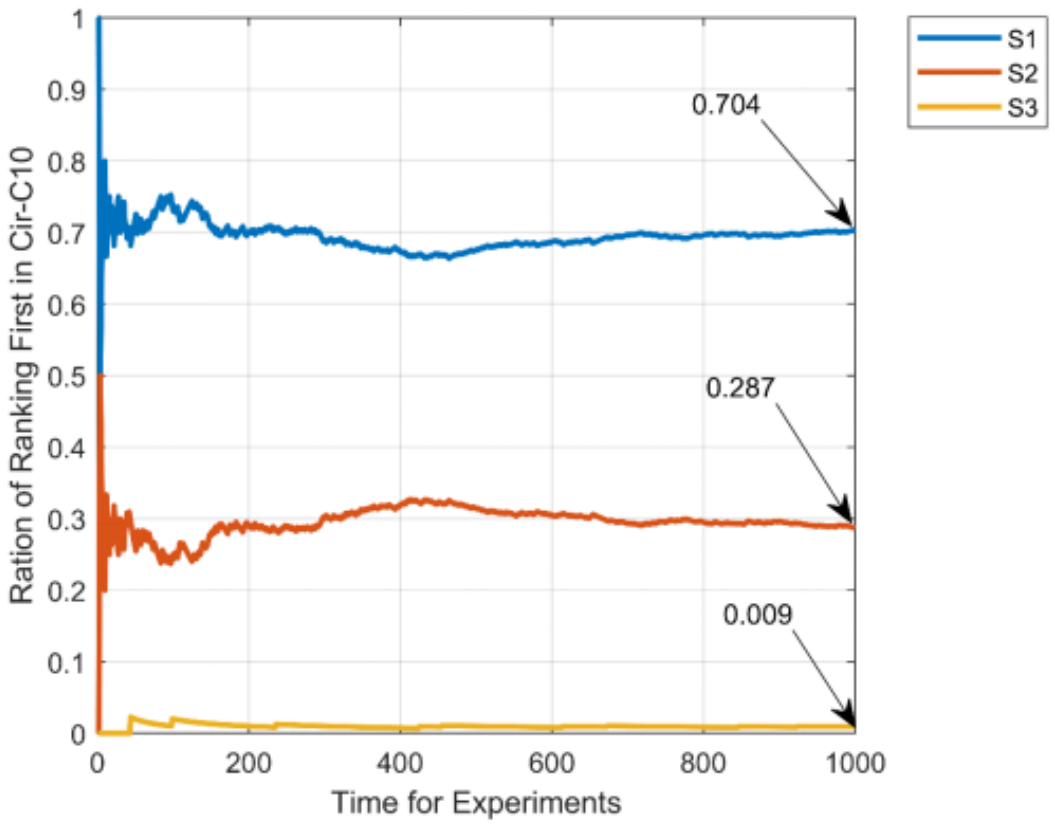

(10) Probability of ranking as first for C10 in circularity

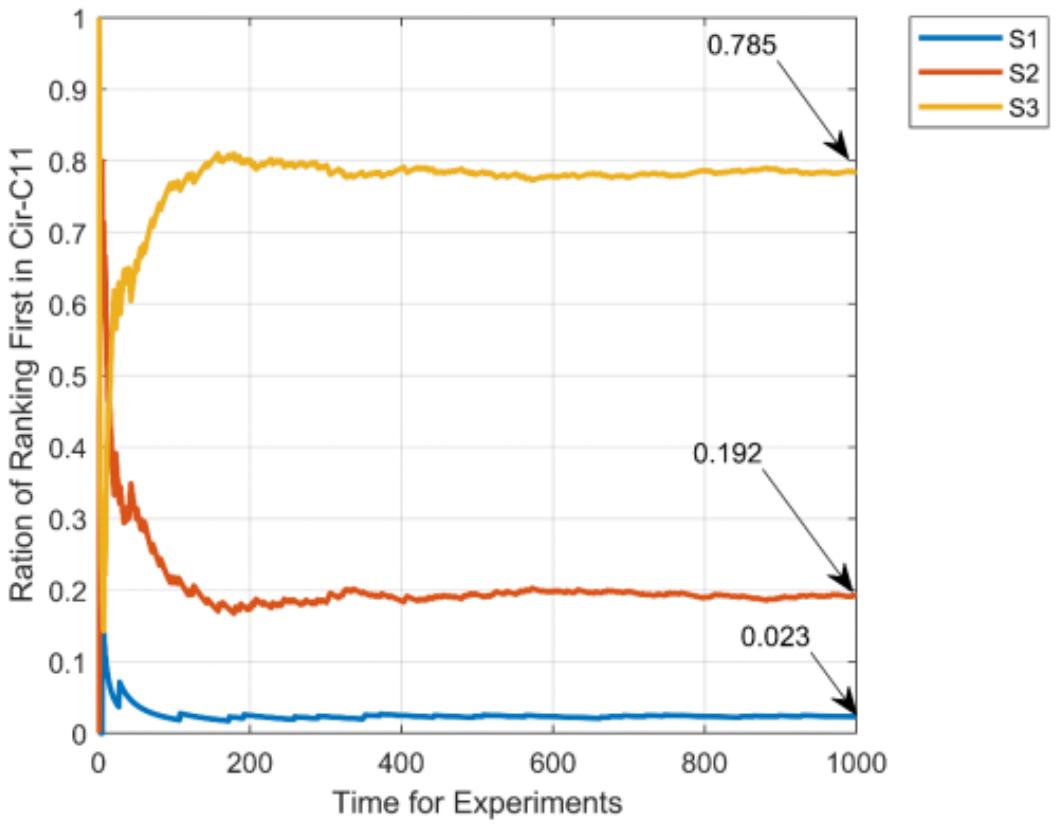

(11) Probability of ranking as first for C11 in circularity

Figure 24 Probability of ranking as first for 12 criteria in Circularity

## 5.4 Analysis

Based on the above probability results, a comparison between the sustainability and circularity analysis for three scenarios is shown in Figure 25.The circularity value is higher than sustainability level for three scenarios.  S3 ranks as first in both sustainability and circularity while S1 ranks as third in both context.

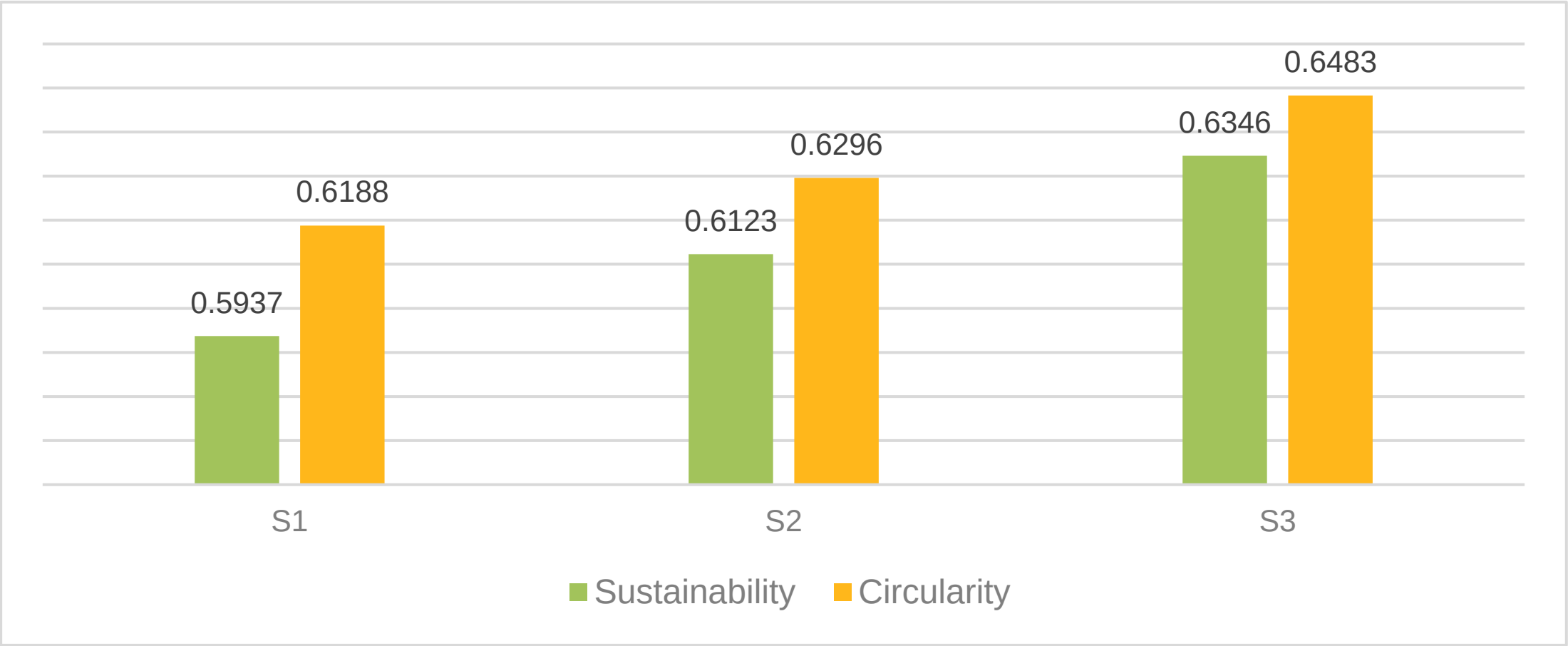

Figure 25 Comparisons for sustainability and circularity of three scenarios

This chart (Figure 26) compares the performance of sustainability and circularity in four aspects. The three scenarios demonstrate a similar trend. The social and technological performances in both sustainability and circularity are comparable. However, circularity has a better environmental performance than sustainability, while sustainability has better economic performance than circularity.

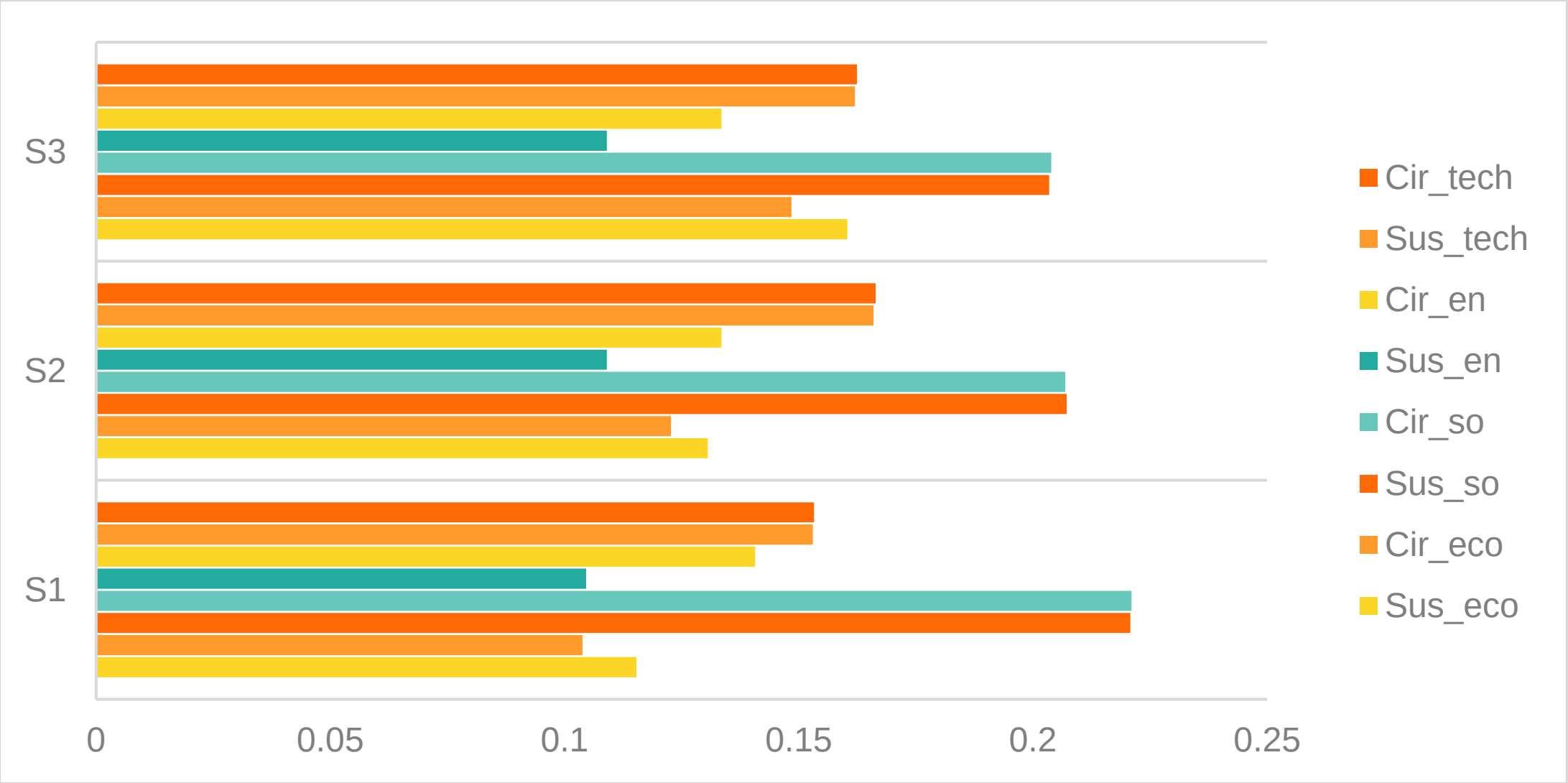

Figure 26 Four requirement level comparisons for three Scenarios

In order to compare the performance of 12 criteria across three scenarios, Sankey diagrams are used to track the performance of each criteria in terms of sustainability and circularity. When analyzing S1, it was found that C3, C4, C6, C7, and C10 had better sustainability performance. In the case of S2, C4, C6, C7, C9, and C10 had better circular performance than sustainability. And in S3, C4, C6, C7, C9, and C10 had better circularity performance than sustainability.

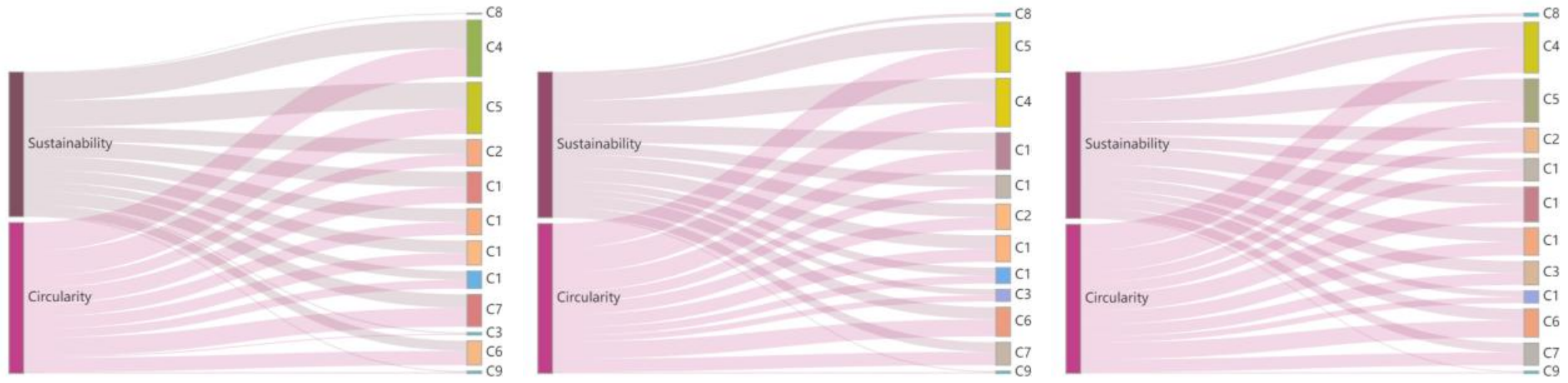

(a) S1 (b) S2 (c) S3

* Waste generation (C8) in sustainability is deleted in circularity, the value of C8 within circularity is set as 0. And the values of C9-C12 are the values of C8-C11

Figure 27 Sankey Diagram for three scenarios

The heatmap (Figure 28) assists in decision-making by providing the likelihood of ranking first in sustainability, circularity, and four other requirements. Among the three options, S3 has the highest probability of ranking first in five selection categories, while S2 has the highest probability in three and S1 in two.

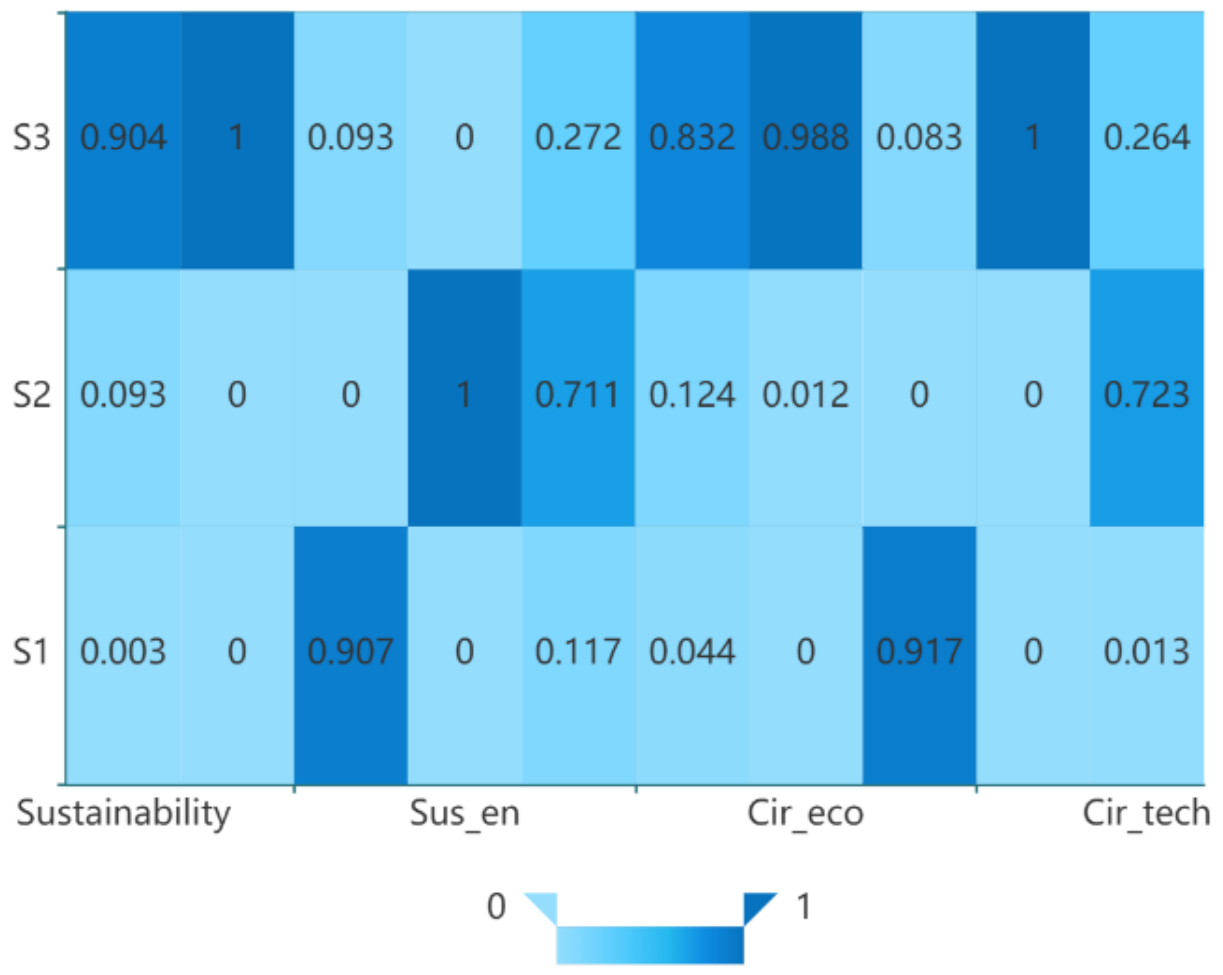

Figure 28 Heatmap for fanking as first for three scenarios

The figure displays a comparison of probabilities for ranking first in 12 sustainability criteria and 11 circularity criteria. The colors range from red to yellow to green, indicating decreasing values. S1 ranks first in six sustainability and circularity criteria, although different ones. S2 has the highest probability in three sustainability and two circularity criteria. S3 performs best in two sustainability and circularity criteria, namely C3 and C12 (C11 for circularity).

| | Sus_C1 | Sus_C2 | Sus_C3 | Sus_C4 | Sus_C5 | Sus_C6 | Sus_C7 | Sus_C8 | Sus_C9 | Sus_C10 | Sus_C11 | Sus_C12 |
|---|---|---|---|---|---|---|---|---|---|---|---|---|
| S1 | 0.226 | 0.582 | 0 | 1 | 0.769 | 0 | 1 | 0 | 1 | 0 | 0.704 | 0.023 |
| S2 | 0.246 | 0.018 | 0 | 0 | 0 | 1 | 0 | 1 | 0 | 1 | 0.287 | 0.192 |
| S3 | 0.528 | 0.4 | 1 | 0 | 0.231 | 0 | 0 | 0 | 0 | 0 | 0.009 | 0.785 |

Figure 29 Heatmap for probability of ranking as first for 12 criteria of sustainability

| | Cir_C1 | Cir_C2 | Cir_C3 | Cir_C4 | Cir_C5 | Cir_C6 | Cir_C7 | Cir_C8 | Cir_C9 | Cir_C10 | Cir_C11 |
|---|---|---|---|---|---|---|---|---|---|---|---|
| S1 | 0.67 | 0.316 | 0 | 1 | 0.769 | 0 | 1 | 1 | 0 | 0.704 | 0.023 |
| S2 | 0 | 0.332 | 0 | 0 | 0 | 1 | 0 | 0 | 1 | 0.287 | 0.192 |
| S3 | 0.33 | 0.352 | 1 | 0 | 0.231 | 0 | 0 | 0 | 0 | 0.009 | 0.785 |

Figure 30 Heatmap for probability of ranking as first for 11 criteria of circularity

## 6. Discussion

AHP's systematic approach derives weights for each criterion on the third layer, encapsulating expert opinion on their relative importance. The homogeneity in the weights across stakeholders lends credibility to the framework and underlines its inclusive, consensus-driven approach. The Monte Carlo simulation, augmented with Latin Hypercube Sampling, is key for uncertainty analysis. By generating numerous combinations of weights for the fourth-layer indicators, it accommodates inherent uncertainties in the assessment process. This range of potential outcomes for the sustainability and circularity performance of the three MMC product scenarios is vital in such complex, uncertain systems. Together, AHP and Monte Carlo simulation bolster the robustness, transparency, and adaptability of the framework, enabling it to handle varied scenarios and adapt to new information. The relative importance of different criteria and the key drivers of uncertainty gleaned from the AHP and Monte Carlo results provide a valuable direction for future research and policy-making within the MMC sector.

Despite the substantial insights derived from the AHP and Monte Carlo simulation methodologies, it's crucial to acknowledge their limitations. AHP's dependence on expert judgment may introduce subjectivity, potentially skewing weight determination. Additionally, its assumption of criteria and indicator independence may not always apply in interrelated systems like the MMC sector. As for Monte Carlo simulation, its reliance on quality input data can limit its effectiveness. The simulation might yield unlikely or unrealistic weight combinations, given its random generation approach.

7. Conclusionprobabilistic method is utilised to generate random weights for the indicator level. By incorporating probabilistic techniques, we can enhance the robustness and objectiveness of the evaluation process, enabling a more reliable and comprehensive understanding of the circularity and sustainability performance in the construction sector. This approach will ultimately contribute to better decision-making and the successful implementation of sustainable and circular practices in the industry.

Based on the results of this paper, other MCDM methods can be applied with Monte Carlo to compare the speed of result convergence to operate sensitivity analysis. More scenarios of MMC can be compared at different levels like product level, structural level and building level.

## 8. Acknowledgement

This research is funded by EPSRC through the Interdisciplinary Circular Economy Centre for Mineral-Based Construction Materials from the UK Research and Innovation (EPSRC Reference: EP/V011820/1).

Data Availability Statement (DAS): Data will be made available on request.